\documentclass[francais, 11pt]{smfart}
\usepackage{smfenum,smfthm}
\usepackage[francais, english]{babel}
\usepackage{amsmath}
\usepackage{amssymb}
\usepackage{amsfonts}
\usepackage{mathrsfs}
\usepackage{graphicx}
\usepackage{epsfig}
\DeclareMathAlphabet{\mathpzc}{OT1}{pzc}{m}{it}

\newtheorem{theoreme}{Th\`eor\`eme}

\newcommand{\partn}[1]{{\smallskip \noindent \textbf{#1.}}}

\renewcommand\P{\mathbb{P}}

\newcommand\R{\mathbb{R}}

\newcommand\cB{\mathcal{B}}

\newcommand\cO{\mathcal{O}}

\newcommand\cS{\mathcal{S}}

\newcommand\fB{\mathfrak{B}}
\newcommand\fC{\mathfrak{C}}
\newcommand\fD{\mathfrak{D}}
\newcommand\fm{\mathfrak{m}}
\newcommand\fF{\mathfrak{F}}
\newcommand\fG{\mathfrak{G}}
\newcommand\fP{\mathfrak{P}}
\newcommand\fS{\mathfrak{S}}
\newcommand\fU{\mathfrak{U}}

\newcommand\sA{\mathscr{A}}
\newcommand\sB{\mathscr{B}}

\newcommand\sP{\mathscr{P}}

\newcommand\hB{\widehat{B}}
\newcommand\hC{\widehat{C}}

\newcommand\hX{\widehat{X}}
\newcommand\hY{\widehat{Y}}

\newcommand\hV{\widehat{V}}
\newcommand\hW{\widehat{W}}

\newcommand\tK{\widetilde{K}}

\newcommand\tY{\widetilde{Y}}

\newcommand\hfB{\widehat{\fB}}
\newcommand\hfF{\widehat{\fF}}

\newcommand\tcS{\widetilde{\cS}}

\newcommand\ov{\overline}

\renewcommand{\=}{ : = }

\DeclareMathOperator{\diam}{diam}

\DeclareMathOperator{\dist}{dist}

\newcommand{\case}[2]{{\smallskip \noindent \textbf{#1.} \textit{\textsf{#2}}}}
\newcommand{\axiomeB}[1]{{\rm ($B_{#1}$)}}

\def\pK{{\P^1_K}}
\def\OK{{\cO_K}}
\def\mK{{\fm_K}}

\def\ptK{{\P^1_{\tK}}}

\def\ber{\mathsf{P}}
\def\berK{{\ber_K^1}}
\def\berV{\mathsf{V}}

\def\arb{{\sA}}
\def\arbX{{\arb_X}}
\def\arbK{{\arb_K}}
\def\arbKR{{\arb_K^\R}}
\def\arbXR{{\arb_X^\R}}
\def\arbp{{\sP}}
\def\arbpK{{\arbp_K}}
\def\arbpKR{{\arbp_K^\R}}
\def\arbboule{{\sB}}
\def\arbbouleR{\arbboule^\R}
\def\distX{{\dist}}
\def\distAX{{\delta}}

\author[J. Rivera-Letelier]{Juan Rivera-Letelier$^\dag$}
\address{Departamento de Matem{\'a}ticas, Universidad Cat{\'o}lica del Norte, Casilla~1280, Antofagasta, Chile}
\email{rivera-letelier@ucn.cl}
\title[Droite projective de Berkovich]{Notes sur la droite projective \\ de Berkovich}

\thanks{$\dag$~Partiellement soutenu par les projets Fondecyt \# 1040683, MeceSup UCN\nobreakdash-0202, South American Network on Dynamical Systems (PROSUL), Nucleo Milenio P01\nobreakdash-005 et MeceSup PUC\nobreakdash-UCH\nobreakdash-0103.}

\begin{document}
\frontmatter


\maketitle
\setcounter{tocdepth}{1}
\tableofcontents

\mainmatter
\section{Introduction}
\subsection{Espaces de Berkovich}
Vers la fin des ann{\'e}es 1980, V.G.~Berkovich a introduit des espaces
analytiques non archim{\'e}diens~\cite{Berrouge, Berblue}, dont la
topologie est plus agr{\'e}able que les espaces analytiques rigides
introduits pr{\'e}c{\'e}demment par J.~Tate~\cite{Ta}.  Contrairement aux
espaces analytiques rigides, les espaces de Berkovich sont connexes
par arcs, localement compacts et localement contractibles.  Ils ont
d{\'e}j{\`a} trouv\'e des applications en g{\'e}om{\'e}trie alg{\'e}brique, physique
math{\'e}matique et en syst{\`e}mes dynamiques.  On pourra consulter~\cite{D}
pour plus de pr{\'e}cisions et de r{\'e}f{\'e}rences.

L'objectif de ces notes est de donner une description {\'e}l{\'e}mentaire de
l'un des espaces de Berkovich le plus simple~: celui associ{\'e} {\`a} la
droite projective, o{\`u} le corps ultram{\'e}trique de base sera suppos{\'e}
complet et alg{\'e}briquement clos.  Cet espace, qu'on appelera par la
suite \textit{droite projective de Berkovich}, poss{\`e}de une structure
d'arbre, o{\`u} chaque sommet poss{\`e}de un nombre infini de branches et o{\`u}
l'ensemble de sommets est dense.

En dehors de son int{\'e}r{\^e}t intrins{\`e}que, la droite projective de Berkovich est int{\'e}ressant car ses parties ouvertes sont des domaines naturels des fonctions \og{}analytiques\fg{}.
Les bonnes propri{\'e}t{\'e}s topologiques de cet espace font que la th{\'e}orie des fonctions analytiques qui en r\'esulte est tr{\`e}s naturelle et proche de la th{\'e}orie complexe, voir~\cite{BR1, Be Ahlfors, Be Ahlfors II, Cherry, CherryRu, injective, elements}.
La droite projective de Berkovich permet aussi d'avoir une th\'eorie du potentiel tout \`a fait analogue \`a la th\'eorie complexe usuelle.
Mentionnons ici~\cite{BR1,FJ}, et notamment~\cite{Th} pour un d\'eveloppement achev\'e dans le cadre plus g\'en\'eral des courbes lisses.

Les espaces de Berkovich sont aussi des espaces naturels pour {\'e}tudier des syst{\`e}mes dynamiques d'origine alg{\'e}brique.
Leur topologie localement compacte garantit des bonnes propri{\'e}t{\'e}s de la th{\'e}orie de la mesure, pour {\'e}tudier des propri{\'e}t{\'e}s ergodiques~\cite{BR2, CL, FR1, FR2}.
Dans le cas particulier de la droite projective, on~a une th{\'e}orie des ensembles Julia et Fatou plus naturelle et plus proche de la th{\'e}orie classique complexe~\cite{Julia/Fatou}, que son analogue sur la droite projective usuelle~\cite{Bez, Hs, these}.
\subsection{Sur ces notes}
On d{\'e}crit ici la droite projective de deux fa{\c c}ons, distinctes de la
d{\'e}finition originale de Berkovich.  La d{\'e}finition de Berkovich est
bas{\'e}e sur des semi-normes multiplicatives d'alg\`ebres de Banach ultram\'etriques, alors que les d{\'e}scriptions
qu'on fait ici sont de nature plus g{\'e}om{\'e}triques.

Premi{\`e}rement on d{\'e}crit la droite projective comme la compl{\'e}tion de la
droite projective usuelle par rapport {\`a} une certaine structure
uniforme.  Cette structure uniforme est diff\'erente de la structure
uniforme induite par la distance chordale, mais elle induit la m{\^e}me
topologie sur la droite projective usuelle.  Il est tr{\`e}s facile de
voir que l'espace qui en resulte est compact.  Cette structure uniforme
est tr{\`e}s utile, car elle est d{\'e}finie m{\^e}me dans le cas o{\`u} la droite
projective de Berkovich n'est pas m{\'e}trisable.  On d{\'e}crit cette
construction dans un contexte plus g{\'e}n{\'e}ral
en~\S~\ref{s:compactification}.  En~\S~\ref{s:compactification droite}
on applique cette construction au cas de la droite projective et on
relie l'espace qui en resulte avec la droite projective de Berkovich.
En~\S~\ref{s:geometrie de berK} on {\'e}tudie quelques propri{\'e}t{\'e}s
g{\'e}om{\'e}triques de cet espace.

Deuxi{\`e}ment on obtient la droite projective de Berkovich, muni de
la \og{}topologie fine\fg{}~, par un proc{\'e}d{\'e} g{\'e}n{\'e}ral et bien connu, qui {\`a} chaque espace ultram{\'e}trique associe un arbre.  La droite
projective de Berkovich c'est l'arbre que l'on obtient par ce proc{\'e}d{\'e},
lorsque l'espace ultram{\'e}trique est {\'e}gal {\`a} la droite projective usuelle
muni de la distance chordale.  La partie \og{}non singuli\`ere\fg{} de
l'espace de Berkovich s'identifie ainsi avec \og{}l'espace des
boules\fg{}.  En~\S~\ref{s:l'arbre d'un espace} on rappelle ce proc{\'e}d{\'e}
g{\'e}n{\'e}ral et en~\S~\ref{s:l'arbre d'un corps} on l'applique au cas
particulier d'un corps ultram{\'e}trique.

Finalement en~\S~\ref{s:topologie fine} on fait le lien entre les deux
constructions.

\subsection{Remarques et r{\'e}f{\'e}rences}
On pourra consulter~\cite{BerICM,D,FvdP,L-S} pour une introduction aux espaces de Berkovich en g\'en\'eral, et~\cite{BR1} pour le cas particulier de la droite projective \'etudi\'ee ici. 
Une th{\'e}orie proche de celle de Berkovich a {\'e}t{\'e} d{\'e}velopp{\'e}e par R.~Huber dans~\cite{Hub}.
Plusieurs \'el\'ements de la th\'eorie de Berkovich ont \'et\'e \'etudi\'es avant les travaux de Berkvoich.
Notamment, M.~van~der~Put a introduit dans~\cite{vdP} des \og{}points g\'eom\'etriques\fg{} d'un espace analytique rigide, appel\'es \og{}filtres premiers\fg{}.
Il a \'et\'e montr\'e apr\`es~\cite{PS} que ces points sont \'etroitement reli\'es \`a l'espace de Berkovich correspondant, voir aussi~\cite{FvdP}.
Les espaces de Berkovich sont model\'es par les espaces de semi-normes multiplicatives d'alg\`ebres de Banach ultram\'etriques.
Ces derniers espaces ont \'et\'e \'etudi\'es notamment par B.~Guennebaud, G.~Garandel et A.~Escassut.
Voir par exemple~\cite{Es80,Ga,Gu}, et les livres plus r\'ecents~\cite{Es95,Es03}.

L'approche aux espaces de Berkovich par les structures uniformes est
apparemment nouvelle.  Cependant ces {\'e}l{\'e}ments sont d{\'e}j{\`a} dans la
litt{\'e}rature d{\'e}puis longtemps~: les filtres de Cauchy minimaux (les
points de la droite projective de Berkovich par cet approche, voir~Proposition~\ref{p:points de berK}) sont exactement les \og{}filtres circulaires\fg{} introduits par Garandel
dans~\cite{Ga} pour caract\'eriser certaines semi-normes
multiplicatives, voir aussi~\cite{Es95,Es03}.
Ils sont aussi proches des \og{}filtres premiers\fg{}, mentionn\'es plus haut.

Diff{\'e}rents aspects topologiques et g\'eom\'etriques de (parties de) la droite projective de Berkovich sont {\'e}tudi\'es dans~\cite{BHM, Mai1, Mai2, hyp, elements}, voir aussi le livre~\cite{Es03}.

Comme il a \'et\'e remarqu\'e par Berkovich dans~\cite[\S~$5$]{Berrouge}, dans certaines cas ces espaces sont \'etroitement reli\'es aux immeubles de Bruhat-Tits, voir aussi~\cite{W}.
En particulier la droite projective de Berkovich est tr\`es reli\'ee \`a l'immeuble de Bruhat-Tits associ\'e au groupe $SL_2$ du corps de base, voir par exemple~\cite[\S\S~$7.2$, $7.3$]{hyp}.
Notons finalement l'\og{}arbre de valuation\fg{} \'etudi\'e dans~\cite{FJ} est assez reli\'e \`a la droite projective de Berkovich.

\subsection{Remerciements}
Je remercie M.~Baker et R.~Rumely pour leur int\'er\^et et leur commentaires lors d'un minicours que j'ai donn\'e \`a l'Universit\'e de Georgia en Avril~2005, lesquels ont \'et\'e tr\`es outils pour am\'eliorer l'exposition de ce travail.
Je remercie aussi M.~Baker, R.~Benedetto et C.~Favre qui ont fait des nombreuses corrections concernant une version pr\'eliminaire de ces notes.
Merci aussi \`a J.~Kiwi pour son int\'er\^et et pour ces nombreux apports lors des discussions qu'on a eu \`a ce sujet.

Ces notes ont \'et\'e \'ecrites au cours des s\'ejours de l'auteur \`a l'Instituto de Matem\'atica Pura e Aplicada (Rio de Janeiro, Brazil), Institute of Mathematics of the Polish Academy of Sciences (Warsaw, Poland), Institut de Math\'ematiques de Jussieu (Paris, France), Centro de Modelamiento Matem\'atico et P.~Universidad Cat\'olica (Santiago, Chile).
Je tiens \`a remercier toutes ces institutions pour leur accueil chaleureux. 
\newpage
\section{Pr{\'e}liminaires}
Soit $K$ un corps muni d'une norme ultram{\'e}trique $| \cdot |$ pour
laquelle $K$ est complet.  Alors l'ensemble
$$
|K^*| \= \{ |z| \mid z \in K \} ~,
$$
est un sous-groupe multiplicatif de~$\R$ qu'on appelle {\it groupe
  des valeurs} de $K$.  On d{\'e}signe par $\OK \= \{ z \in K \mid |z|
\le 1 \}$ l'anneau des entiers de $K$ et par $\mK \= \{ z \in K \mid
|z| < 1 \}$ l'id{\'e}al maximal de $\OK$.  On notera $\tK \= \OK / \mK$ le
corps r{\'e}siduel de $K$.

\subsection{La droite projective.}
On note par $\pK$ la droite projective de~$K$, qui est l'ensemble des
droites dans $K \times K$ passant par $(0, 0)$.  Pour $(x, y) \in K
\times K \setminus \{ (0, 0) \}$, on d{\'e}signe par $[x : y] \in \pK$ le
point correspondant {\`a} la droite $\{ (\lambda x, \lambda y) \mid
\lambda \in K \}$.  On d{\'e}signe par $\infty$ le point $[1 : 0] \in
\pK$, et on identifie $K$ {\`a} $\pK \setminus \{ \infty \}$ par
l'application $x \mapsto [x : 1]$.

On {\'e}tend la projection de $\OK$ vers $\tK$ en une projection de $\pK$
vers $\ptK$, de telle fa{\c c}on que l'ensemble $\pK \setminus \OK$ se
projette dans $\infty \in \ptK$.  Les fibres de cette projection seront
appel{\'e}s {\it classes r\'esiduelles}.

La fonction $\Delta : \pK \times \pK \to [0, 1]$ d{\'e}finie par
$$
\Delta([x:y],[x':y']) = \frac{|xy' - yx'|}{\max \{ |x|, |y| \}
  \cdot \max \{ |x'|, |y'| \} },
$$
d{\'e}finit une distance sur $\pK$ qu'on appelle \textit{distance chordale}.
Notons que pour $z, z' \in K \subset \pK$ on~a
$$
\Delta(z, z') = \frac{|z - z'|}{\max \{ 1, |z| \} \cdot \max \{ 1, |z'| \}}.
$$
\subsection{Boules}
Lorsque $B$ est un sous-ensemble de $K$ on pose $\diam(B) = \sup_{z, z' \in B} \{ |z - z'| \}$.

Etant donn{\'e}s $r \in |K^*|$ et $a \in K$, on appelle les ensembles
$$
        \{ z \in K \mid |z - a| < r \} \mbox{ et }
                \{ z \in K \mid |z - a| \le r\}
                $$
                \textit{boule ouverte de} $K$ et \textit{boule
                  ferm{\'e}e de} $K$, respectivement.  Lorsque $r \not \in
                |K^*|$, ces deux ensembles co{\"\i}ncident et constituent
                ce qu'on appelle une \textit{boule irrationnelle de}
                $K$.  Notons que par d{\'e}finition une boule $B$ de $K$
                est irrationnelle si et seulement si $\diam(B) \not
                \in |K^*|$~; en particulier, \textit{si $B$ est
                  ouverte ou ferm{\'e}e, alors $\diam(B) \in |K^*|$}.  Si
                deux boules de $K$ s'intersectent, alors l'une est
                contenue dans l'autre.
                
                Une \textit{boule ouverte} (resp. \textit{ferm{\'e}e,
                  irrationnelle}) de $\pK$ est soit une boule de $K$
                de m{\^e}me nature, soit le compl{\'e}mentaire d'une boule
                ferm{\'e}e (resp. ouverte, resp. irrationnelle) de $K$.
                \textit{Dans ce qui suit le mot boule d{\'e}signera une
                  boule de~$\pK$.}
\subsection{Affino{\"\i}des}\label{ss:affinoides}
Un \textit{affino{\"\i}de ferm{\'e}} (resp. \textit{ouvert}) est une
intersection finie de boules ferm{\'e}es (resp. ouvertes).  L'ensemble
$\pK$ est un affino{\"\i}de ferm{\'e} (resp. ouvert), car il l'est une
intersection vide de boules ferm{\'e}es (resp. ouvertes).  Une
intersection finie non vide d'affino{\"\i}des ferm{\'e}s (resp. ouverts) est un
affino{\"\i}de ferm{\'e} (resp. ouvert).  La r{\'e}union de deux affino{\"\i}des ferm{\'e}s
(resp. ouverts) dont l'intersection est non vide est un affino{\"\i}de
ferm{\'e} (resp. ouvert).
\begin{rema}
La terminologie qu'on utilise ici est distincte \`a celle de la g\'eom\'etrie rigide.
En g\'eom\'etrie rigide un \og{}affino\"{\i}de\fg{} est ce qu'on appelle ici une reunion finie d'affino\"{\i}des ferm\'es, et un \og{}affino\"{\i}de connexe\fg{} correspond \`a un affino\"{\i}de ferm\'e. 
\end{rema}

\begin{prop}\label{p:affinoides}
  Soit $n$ un entier positif et pour chaque $j = 1, \ldots, n$ soit
  $X_j$ un affino{\"\i}de ouvert.  Pour chaque $x \in X = X_1 \cup \cdots
  \cup X_n$ l'union $Y$ de tous les affino{\"\i}des ouverts contenant $x$ et
  contenus dans $X$ est un affino{\"\i}de ouvert.  On appelle $Y$ la {\sf
    composante de $X$ contenant} $x$.  Les composantes de $X$ sont
  disjointes deux {\`a} deux et il n'y a qu'un nombre fini d'eux.
\end{prop}
\begin{proof}
  Soit $J \subset \{ 1, \ldots, n \}$ l'ensemble des $j$ tel qu'il
  existe un affino{\"\i}de ouvert contenu dans $X$, contenant $x$ et qui
  rencontre $X_j$.  Posons $Y = \bigcup_{j \in J} X_j$ et soit $Z$ un affino{\"\i}de
  ouvert contenant $x$ et contenu dans $X$.  Chaque point de $Z
  \subset X$ est contenu dans l'un des $X_j$ o{\`u} $j \in J$ par
  d{\'e}finition de $J$, d'o{\`u} $Z \subset Y$.
  
  D'autre part, pour chaque $j \in J$ il existe un affino{\"\i}de ouvert
  $Z_j$ contenu dans $X$, qui contient $x$ et qui rencontre $X_j$.  Par
  cons{\'e}quent $Z_j \cup X_j$ est un affino{\"\i}de ouvert ayant les m{\^e}mes
  propri{\'e}t{\'e}s, d'o{\`u} $Y \subset \bigcup_{j \in J} ( Z_j \cup X_j ) \subset
  Y$.  Par cons{\'e}quent l'affino{\"\i}de ouvert $Y$ est {\'e}gal {\`a} l'union de
  tous les affino{\"\i}des ouverts contenus dans $X$ et contenant $x$.
  
  Si la composante $Y'$ de $X$ contenant un point $x' \in X$ rencontre
  $x$, alors $Y \cup Y'$ est un affino{\"\i}de ouvert contenu dans $X$ et
  contenant $x$ et $x'$, d'o{\`u} $Y = Y'$.  Finalement, comme chaque
  composante de $X$ est de la forme $\bigcup_{j \in J} X_j$, on conclut que $X$
  ne poss{\`e}de qu'un nombre fini de composantes.
\end{proof}

Le corollaire suivant est imm{\'e}diat.

\begin{coro}\label{c:affinoides}
  Soit $X$ une union finie d'affino{\"\i}des ouverts et soit $Y$ un
  affino{\"\i}de ouvert contenu dans $X$.  Alors $Y$ est contenu dans l'une
  des composantes de $X$.
\end{coro}

\newpage 
\section{Compactification des recouvrements finis}\label{s:compactification}
Apr\`es quelques rappels sur les filtres (\S~\ref{ss:rappel filtres}) et sur les structures uniformes (\S~\ref{ss:rappel structures uniformes}), on d\'ecrit un proc\'ed\'e qui \`a chaque base convenable d'un espace topologique associe une structure uniforme compatible avec la topologie originale (\S\S~\ref{ss:structure uniforme}, \ref{ss:topologie deduite}).
Apr\`es on d\'ecrit le s\'epar\'ee compl\'ete de cet espace uniforme (\S~\ref{ss:completion}) et on \'etudie quelques propri\'et\'es topologiques de cet espace (\S~\ref{ss:topologie du complete}).
\subsection{Rappel sur les filtres}\label{ss:rappel filtres}
Un \textit{filtre} sur un ensemble $X$ est une collection $\fF$ de
parties de $X$ satisfaisant les propri{\'e}t{\'e}s suivantes.
\begin{enumerate}
\item[1.]
Toute partie de $X$ contenant un {\'e}l{\'e}ment de $\fF$ appartient {\`a}~$\fF$.
\item[2.]
Toute intersection finie d'{\'e}l{\'e}ments de $\fF$ appartient {\`a}~$\fF$.
\item[3.]
L'ensemble vide $\emptyset$ n'appartient pas {\`a}~$\fF$.
\end{enumerate}
Notons que les propri{\'e}t{\'e}s~$2$ et~$3$ impliquent que toute intersection
finie d'{\'e}l{\'e}ments de $\fF$ est non vide.
Une collection $\fB$ de parties de $X$ est une \textit{base de filtre}, si $\emptyset \not \in \fB$ et si toute intersection finie d'\'el\'ements de $\fB$ contient l'un des \'el\'ements de $\fB$.
Lorsque $\fB$ est une base d'un filtre, la collection
$$
\{ Y \in X \mid \text{il existe $B \in \fB$ tel que $B \subset Y$} \}
$$
est un filtre sur $X$ qu'on appelle \textit{le filtre sur $X$ engendr\'e par $\fB$}.

Etant donn\'es deux filtres $\fF$ et $\fF'$ sur un ensemble~$X$, on dit
que $\fF$ est \textit{plus fin} que $\fF'$, si $\fF' \subset \fF$.
Un \textit{ultrafiltre} sur un ensemble~$X$ est un filtre tel qu'il
n'existe aucun filtre strictement plus fin que lui.
Il est facile de voir qu'un filtre~$\fF$ est un ultrafiltre si et seulement si pour chaque sous-ensemble~$Y$ de~$X$ on a $Y \in \fF$ ou $X \setminus Y \in \fF$.
Pour tout filtre $\fF$ sur un ensemble $X$ il existe un ultrafiltre moins fin que $\fF$
\cite[I, \S6, Th{\'e}or{\`e}me~$1$]{Bou}.  Il est facile de voir que si $\fF$
est un ultrafiltre sur un ensemble $X$ et si $(Y_i)_{i = 1, \ldots n}$
est un recouvrement fini de~$X$, alors l'un au moins des $Y_i$
appartient {\`a}~$\fF$, voir~\cite[I, \S~6, corollaire de la
Proposition~5]{Bou}.

Lorsque $X$ est muni d'une
topologie, on dit qu'un filtre $\fF$ sur $X$ \textit{converge vers un
  point $x$} de $X$, s'il contient tout voisinage de~$x$.
Dans ce cas on dit que~$x$ \textit{est le point limite de} $\fF$.
\begin{prop}\label{p:compacite via ultrafiltres}
  Un espace topologique s{\'e}par{\'e} $X$ est compact si et seulement si tout
  ultrafiltre sur $X$ est convergent.
\end{prop}
\begin{proof}
  Supposons d'abord qu'il existe un ultrafiltre~$\fF$ sur~$X$ qui
  n'est pas convergent.  Alors la collection $\{ X \setminus \ov{Y}
  \mid Y \in \fF \}$ est un recouvrement ouvert de~$X$ ne poss{\`e}dant
  aucun sous-recouvrement fini.

Supposons d'autre part que tout ultrafiltre sur~$X$ soit convergent et
supposons par l'absurde qu'il existe un recouvrement ouvert $\fG$ de $X$
ne possedant aucun sous-recouvrement fini.  Alors la
collection $\{ X \setminus O \mid O \in \fG\}$ forme une base d'un
filtre~$\fF$ sur~$X$.  Si $\fF'$ est un ultrafiltre sur $X$ moins fin
que~$\fF$, alors $\fF'$ est convergent par hypoth{\`e}se.  Mais le point
limite de $\fF'$ appartient {\`a} tous les elements de~$\fF$, ce qui
contredit le fait que $\fG$ est un recouvrement de~$X$.
\end{proof}
\subsection{Rappel sur les espaces uniformes}\label{ss:rappel structures uniformes}
On note par $\Delta_X \= \{ (x, x) \in X \times X \mid x \in X \}$ la
diagonale dans $X \times X$ et pour deux parties $V$ et $V'$ de $X
\times X$ on pose
$$
V^{-1} \= \{ (x, x') \in X \times X \mid (x', x) \in V \} ~,
$$
\begin{multline*}
V \circ V' \= \{ (x, x') \in X \times X \mid  \\
\text{ il existe $y \in X$ tel que $(x, y) \in V$ et $(y, x') \in V'$} \} ~.
\end{multline*}

Une \textit{structure uniforme} sur un ensemble $X$ est la donn{\'e}e
d'une collection $\fU$ de parties de $X \times X$ satisfaisant les
propri{\'e}t{\'e}s~$1$ et~$2$ des filtres, et les propri{\'e}t{\'e}s suivantes.
\begin{enumerate}
\item[1.]
Tout {\'e}l{\'e}ment de $\fU$ contient la diagonale $\Delta_X$.
\item[2.]
Pour chaque $V \in \fU$ on~a $V^{-1} \in \fU$.
\item[3.]
Pour tout $V \in \fU$ il existe $W \in \fU$ tel que $W \circ W \subset V$.
\end{enumerate}
Les {\'e}l{\'e}ments de $\fU$ sont appel{\'e}s \textit{entourages}.
Notons que par la propri\'et\'e~$1$, la collection $\fU$ est un filtre sur $X \times X$.
On dira qu'une base de filtre $\fB$ sur $X \times X$ est \textit{une base de structure uniforme}, si tout \'el\'ement de $\fB$ contient $\Delta_X$, et si pour tout $V \in \fB$ il existe $V' \in \fB$ tel que $V' \subset V^{-1}$, et $W \in \fB$ tel que $W \circ W \subset V$.
Le filtre engendr\'e par une base de structure uniforme $\fB$ est une structure uniforme, qu'on appelle \textit{la structure uniforme engendr\'ee par $\fB$}.

Une structure uniforme sur un ensemble $X$ induit une topologie sur $X$, caract\'eris{\'e}e par le fait que les voisinages d'un point $x \in X$ sont les ensembles de la forme,
\begin{equation}\label{e-voisinage}
V(x) \= \{ x' \in X \mid (x, x') \in V \} ~,
\end{equation}
o{\`u} $V$ parcourt les entourages de $X$, voir~\cite[II, \S~$1$, Proposition~$1$]{Bou}.
Les ouverts de $X$ sont les ensembles qui sont un voisinage de chacun de ses {\'e}l{\'e}ments.  Lorsque la topologie induite est s{\'e}par{\'e}e, elle est r{\'e}guli{\`e}re, voir~\cite[II, \S~$1$, Proposition~$3$]{Bou}.

Une distance $\dist$ sur $X$ induit une structure uniforme sur $X$, engendr\'ee par la base de structure uniforme constitu\'ee des ensembles
$$
\{ (x, x') \in X \times X \mid \dist(x, x') < r \},
$$
lorsque $r$ parcourt les nombres r\'eels positifs.
La topologie sur $X$ induite par cette structure uniforme co\"{\i}ncide avec celle induite par la distance $\dist$.

Pour les espaces uniformes il y a une notion de compl\'etion qui g\'en\'eralise celle des espaces m\'etriques.
Elle est d\'efinie comme suit.
Lorsque $X$ est un espace uniforme et $V$ est un entourage de~$X$, on dira qu'une partie $Y$ de $X$ est $V$-\textit{petite}, si $Y \times Y \subset V$.
De plus, on dira qu'un filtre est \textit{de Cauchy}, si pour tout entourage $V$ de~$X$ il contient un {\'e}l{\'e}ment $V$-petit, et on dira que l'espce uniforme $(X, \fU)$ est \textit{complet} si tout filtre de Cauchy est convergent.
Pour d\'efinir le \textit{s\'epar\'e compl\'et\'e} $\hX$ d'un espce uniforme~$X$, notons d'abord que pour tout filtre Cauchy il existe un unique filtre de Cauchy moins fin, et qui est minimal avec cette propri\'et\'e~\cite[II, \S~$3$, Proposition~$5$]{Bou}.
Alors $\hX$ est d\'efinit comme l'ensemble de tous les filtres de Cauchy minimaux, muni de la structure uniforme de toutes les parties de $\hX \times \hX$ contenant un ensemble de la forme
\begin{multline}\label{e-entourage de hX}
\hV \= \{ (\fF, \fF') \in X \times X \mid 
\\
\text{ $\fF$ et $\fF'$ contiennent un ensemble $V$-petit commun} \} ~,
\end{multline}
o\`u $V$ est un entourage de~$X$.
Cet espace est complet et s\'epar\'e, et par cons\'equent r\'egulier.
Lorsque l'espace uniforme $X$ est lui-m\^eme s\'epar\'e, il s'identifie au sous-espace partout dense de $\hX$ des filtres convergents~; chaque point de $X$ s'identifie \`a l'unique filtre minimal moins fin que le filtre des voisinages de ce point.
\subsection{Structure uniforme des recouvrements finis}\label{ss:structure uniforme}
Soit $X$ un espace topologique et soit $\fB$ une collection de parties
ouvertes de $X$.  On note par $\fC$ la collection de tous les
recouvrements finis de $X$ par {\'e}l{\'e}ments de $\fB$.  Par la suite on
supposera que $\fB$ satisfait les propri{\'e}t{\'e}s suivantes.
\begin{enumerate}
\item[($B_I$)]\label{BI} Toute intersection finie non vide d'{\'e}l{\'e}ments
  de $\fB$ appartient {\`a}~$\fB$.
\item[($B_{II}$)]\label{BII} Pour tout recouvrement $C$ de $X$ dans
  $\fC$ il existe un recouvrement $C'$ dans $\fC$, satisfaisant la
  propri{\'e}t{\'e} suivante.  Pour chaque point $x$ de $X$ la r{\'e}union de tous
  les {\'e}l{\'e}ments de $C'$ contenant~$x$ est contenue dans l'un des
  {\'e}l{\'e}ments de~$C$.
\end{enumerate}

Pour un recouvrement $C$ dans $\fC$ on pose,
\begin{equation}
V(C) \= \bigcup_{Y \in C} Y \times Y ~.
\end{equation}
\begin{prop}\label{p:structure uniforme}
  La collection $\fU$ d\'efinie par,
$$
\fU = \{ V \subset X \times X \mid \text{il existe $C \in \fC$ tel que $V(C) \subset V$} \},
$$
d{\'e}finit une structure uniforme sur $X$.
\end{prop}
\begin{proof}
  Montrons d'abord que $\fU$ est invariante par l'intersection finie.
  Etant donn{\'e}s des recouvrements $C_1, \ldots, C_n$ dans $\fC$, la
  propri{\'e}t{\'e} \axiomeB{I} implique que le recouvrement fini de $X$
  d{\'e}fini par
$$
C \=
\{ Y_1 \cap \ldots \cap Y_n \neq \emptyset
\mid Y_j \in C_j, j = 1, \dots, n \} ~,
$$
appartient {\`a}~$\fC$.  Il est facile de voir qu'on a $V(C) \subset
V(C_1) \cap \ldots \cap V(C_n)$ et par cons{\'e}quent ce dernier ensemble
appartient {\`a}~$\fU$.

D'autre part, notons que tout {\'e}l{\'e}ment de $\fU$ contient la diagonale
et comme pour chaque $C \in \fC$ l'ensemble $V(C)$ est invariant par
l'involution $(x, x') \mapsto (x', x)$, pour tout $V \in \fU$ on~a
$V^{-1} \in \fU$.  Pour v{\'e}rifier que $\fU$ est une structure uniforme,
il reste {\`a} montrer que pour tout $V \in \fU$ il existe $W \in \fU$ tel
que $W \circ W \subset V$.  Etant donn{\'e} $V \in \fU$, soit $C \in \fC$
tel que $V(C) \subset V$ et soit $C'$ le recouvrement dans $\fC$ donn{\'e}
par la propri{\'e}t{\'e}~\axiomeB{II}.  Etant donn{\'e}s $x, y, z \in \pK$ tels
que $(x, y), (y, z) \in V(C')$, soient $\tY, \tY' \in C'$ tels que $x,
y \in \tY$ et $y, z \in \tY'$.  Si $Y \in C$ contient $\tY$ et $\tY'$,
alors on~a $(x, z) \in Y \times Y \subset V(C)$.
\end{proof}
\subsection{Topologie d{\'e}duite}\label{ss:topologie deduite}
Notons que pour $C \in \fC$ et $x \in X$, l'ensemble $V(C)(x)$ d{\'e}finit
au~\eqref{e-voisinage} pour l'entourage $V \= V(C)$, est {\'e}gal {\`a} la
r{\'e}union des {\'e}l{\'e}ments de $C$ contenant~$x$.

Notons pour r{\'e}f{\'e}rence que lorsque $\fB$ est d{\'e}nombrable, $\fC$ est
d{\'e}nombrable et donc la topologie sur $X$ d{\'e}finie par la structure
uniforme $\fU$ est m{\'e}trisable \cite[IX, \S~$4$, Th{\'e}or{\`e}me~$1$]{Bou}.
\begin{prop}\label{p:topologie deduite}
  La topologie sur~$X$ d{\'e}duite de la structure uniforme sur~$X$
  d{\'e}finie par~$\fU$, est moins fine que la topologie originale de~$X$.
  Ces topologies co{\"\i}ncident, si et seulement si $\fB$ satisfait la
  propri{\'e}t{\'e} suivante.
\begin{quote}
\axiomeB{III} Pour chaque $x \in X$ et chaque voisinage $Y$ de
  $x$ pour la topologie originale de $X$, il existe un recouvrement
  $C$ dans $\fC$ tel que la r{\'e}union des {\'e}l{\'e}ments de $C$ contenant $x$
  est contenue dans $Y$.
\end{quote}
\end{prop}
\begin{proof}
  Pour montrer la premi{\`e}re assertion, notons que pour un recouvrement
  $C$ dans $\fC$ et pour $x \in X$, l'ensemble $V(C)(x) = \bigcup_{Y \in
    C, \; x \in Y} Y$ est une partie ouverte de $X$.
  
  Supposons maintenant $\fB$ satisfait la propri{\'e}t{\'e} \axiomeB{III} et
  soient $x \in X$ et $Y$ un voisinage de~$x$ pour la topologie
  originale.  Il existe alors un recouvrement $C$ dans $\fC$ tel que
  la r{\'e}union des {\'e}l{\'e}ments dans $C$ contenant $x$ soit contenue dans
  $Y$.  C'est {\`a} dire qu'on a $V(C)(x) \subset Y$ et par cons{\'e}quent $Y$
  est aussi un voisinage de $x$ pour la topologie d{\'e}duite de la
  structure uniforme.
  
  Supposons d'autre part que la topologie d{\'e}duite de la structure
  uniforme co{\"\i}ncide avec la topologie originale.  Etant donn{\'e}s $x \in
  X$ et un voisinage $Y$ de $x$, soit $V$ un entourage de $X$ tel
  qu'on~ait $V(x) = Y$.  Soit de plus $C$ un recouvrement dans $\fC$
  tel que $V(C) \subset V$.  Alors la r{\'e}union $V(C)(x)$ des {\'e}l{\'e}ments
  de $C$ contenant $x$ est contenue dans $V(x) = Y$.
\end{proof} 
\subsection{Compl{\'e}tion}\label{ss:completion}
Soit $X$ un espace topologique muni de la structure uniforme d\'efinie \`a partir d'une collection de parties ouvertes satisfaisant les propri{\'e}t{\'e}s \axiomeB{I} et \axiomeB{II}.

On consid\`ere d'abord le lemme suivant.
\begin{lemm}\label{l:filtres de Cauchy}
  Pour qu'un filtre soit de Cauchy il faut et il suffit qu'il contienne
  au moins un {\'e}l{\'e}ment de chaque recouvrement dans $\fC$.
\end{lemm}
\begin{proof}
  Pour voir que la condition est suffisante, notons simplement que si
  $V$ est un entourage et si $C$ est un recouvrement dans $\fC$ tel
  que $V(C) \subset V$, alors tout {\'e}l{\'e}ment de $C$ est $V$-petit.

Pour montrer que la condition est n\'ecessaire, soit $\fF$ un filtre
de Cauchy et soit $C$ un recouvrement dans $\fC$.  Soit $C' \in \fC$
le recouvrement donn{\'e} par la propri{\'e}t{\'e} \axiomeB{II}.  Comme $\fF$ est
un filtre de Cauchy, il existe un ensemble $Y_0 \in \fF$ tel que $Y_0
\times Y_0 \subset V(C')$.  Si l'on fixe $y_0 \in Y_0$, alors la
r{\'e}union $Y$ de tous les {\'e}l{\'e}ments de $C'$ contenant $y_0$, contient
$Y_0$.  Mais $Y$ est contenu dans un certain {\'e}l{\'e}ment de $C$.  Cet
{\'e}l{\'e}ment de $C$ appartient donc {\`a}~$\fF$.
\end{proof}

On d{\'e}signe par $\hX$ le s{\'e}par{\'e} compl{\'e}t{\'e} de $X$.
Etant donn{\'e}e une partie $Y$ de~$X$ on pose
$$
\hY \= \{ \fF \in \hX \mid Y \in \fF \} ~.
$$
Le Lemme~\ref{l:filtres de Cauchy} implique que pour chaque
recouvrement $C \in \fC$, la collection $\hC \= \{ \hY \mid Y \in C
\}$ est un recouvrement de $\hX$.
Notons d'autre part que l'entourage $\hV$ de $\hX$, d\'efini en~\eqref{e-entourage de hX} lorsque $V \= V(C)$, est {\'e}gal {\`a}
$$
\hV(C) \= \bigcup_{Y \in C} \hY \times \hY ~.
$$
Comme les ensembles $V(C)$, lorsque $C$ parcourt $\fC$, constituent une
base de la structure uniforme de $X$, on conclut que les ensembles
$\hV(C)$, lorsque $C$ parcourt $\fC$, constituent une base de la structure
uniforme de $\hX$.

\subsection{Topologie du compl\'et\'e}\label{ss:topologie du complete}
On reprend les notations du paragraphe pr\'ec\'edent.
Apr\`es avoir montr\'e que $\hX$ est compact (Proposition~\ref{p:le complete est compact}), on decrira un base de topologie de cet espace (Proposition~\ref{p:base de la topologie}).

\begin{prop}\label{p:le complete est compact}
L'espace uniforme $\hX$ est compact.
\end{prop}
\begin{proof}
  D'apr{\`e}s la Proposition~\ref{p:compacite via ultrafiltres} il suffit
  de montrer que tout ultrafiltre sur $\hX$ est convergent.  Comme
  $\hX$ est complet, il suffit de montrer que tout ultrafiltre sur
  $\hX$ est de Cauchy.  Soit alors $\hfF$ un ultrafiltre sur $\hX$ et
  soit $\hV$ un entourage de $\hX$.  Soit $C$ un recouvrement dans
  $\fC$ tel que $\bigcup_{Y \in C} \hY \times \hY \subset \hV$.  Comme
  $\hC = \{ \hY \mid Y \in C \}$ est un recouvrement fini de $\hX$,
  l'ultrafiltre $\hfF$ contient l'un des {\'e}l{\'e}ments de $\hC$,
  voir~\S\ref{ss:rappel filtres}.  Comme tous les {\'e}l{\'e}ments de $\hC$
  sont $\hV$-petits, on conclut que $\hfF$ est de Cauchy.
\end{proof}

Soit $Y$ une partie quelconque de $X$.
Alors il est facile de voir que la fermeture topologique $\ov{Y}$ de $Y$ dans $\hX$ est {\'e}gale {\`a}
$$
\ov{Y}
=
\{ \fF \in \hX \mid \text{tout {\'e}l{\'e}ment de $\fF$ rencontre $Y$} \} ~.
$$
Notons qu'on a $Y \subset \hY \subset \ov{Y}$, et par cons{\'e}quent la
fermeture topologique de $\hY$ est {\'e}gale {\`a} $\ov{Y}$.
D'autre part, si $Y$ et $Y'$ sont des parties de $X$, alors il est
facile de voir qu'on~a
\begin{equation}\label{e-core}
\widehat{Y \cap Y'} = \hY \cap \widehat{Y'}
\ \text{ et } \
\hY \cup \widehat{Y'} \subset \widehat{Y \cup Y'} ~.
\end{equation}
En particulier, lorsque les ensembles $Y$ et $Y'$ sont disjoints, les
ensembles $\hY$ et $\widehat{Y'}$ sont aussi disjoints.
\begin{prop}\label{p:base de la topologie}
\

\begin{enumerate}
\item[1.]
  Lorsque $Y$ et $Y'$ sont des parties compl{\'e}mentaires de $X$, les
  ensembles $\hY$ et $\ov{Y'}$ sont des parties compl{\'e}mentaires de $\hX$.
  En particulier l'ensemble $\hY$ est une partie ouverte de $\hX$.
\item[2.]
La collection des parties ouvertes de $\hX$,
$$
\hfB \= \{ \hY \mid Y \in \fB \} ~,
$$
est une base de la topologie de $\hX$.
\end{enumerate}
\end{prop}
\begin{proof}
\

\partn{1}
  Si $\fF \in \hY$ alors $Y$ est un {\'e}l{\'e}ment de $\fF$ disjoint de $Y'$
  et donc $\fF \not \in \ov{Y'}$.  D'autre part, $\fF \not \in
  \ov{Y'}$ implique qu'il existe $Z \in \fF$ disjoint de $Y'$.  On a
  alors $Z \subset Y$, $Y \in \fF$ et $\fF \in \hY$.

\partn{2}
  Soit $\fF$ un point de $\hX$ et soit $\hW$ un voisinage de $\fF$
  dans $\hX$.  Soit $\hV$ un entourage de $\hX$ tel que $\hW =
  \hV(\fF)$ et soit $C$ un recouvrement dans $\fC$ tel que $\hV(C)
  \subset \hV$.  Comme $\{ \hY \mid Y \in C \}$ est un recouvrement de
  $\hX$, il existe $Y \in C$ tel que $\fF \in \hY$.  Alors $\hY$ est
  un {\'e}l{\'e}ment de $\hfB$ contenant $\fF$ et contenu dans $\hV(C)(\fF)
  \subset \hW$.
\end{proof}

\newpage
\section{Compactification de la droite projective}\label{s:compactification droite}
Pour toute cette section on fixe un corps $K$ muni d'une norme
ultram{\'e}trique $| \cdot |$ pour laquelle $K$ est complet.  On suppose
de plus que le corps r{\'e}siduel de $K$ est infini et que le groupe de
valuation de $K$ n'est pas discret.  Ces derni{\`e}res conditions sont
satisfaites, par exemple, lorsque $K$ est alg{\'e}briquement clos et la
norme $| \cdot |$ est non triviale.

Dans cette section on appliquera la construction g{\'e}n{\'e}rale d{\'e}crite
en~\S~\ref{s:compactification}, au cas o{\`u} l'ensemble $X$ est {\'e}gal {\`a}
$\pK$ et $\fB$ est {\'e}gal {\`a} la collection des affino{\"\i}des ouverts de
$\pK$ (\S\S~\ref{ss:collection des affinoides}, \ref{ss:berK}).  De
plus, on d{\'e}crira les points du s{\'e}par{\'e} compl{\'e}t{\'e} correspondant
(\S\S~\ref{ss:points de berK}, \ref{ss:caracteristique minimalite}) et
on reliera cet espace {\`a} la droite projective de Berkovich (\S~\ref{ss:lien}).
\subsection{La collection des affino{\"\i}des ouverts}\label{ss:collection des affinoides}
Dans ce paragraphe on montre que la collection des affino{\"\i}des ouverts
de $\pK$ v{\'e}rifie les propri{\'e}t{\'e}s \axiomeB{I}, \axiomeB{II} et
\axiomeB{III} d{\'e}crites en~\S\S~\ref{ss:structure uniforme},
\ref{ss:topologie deduite}.

Rappelons d'abord que chaque affino{\"\i}de ouvert est une partie ouverte
de $\pK$.  Comme toute intersection finie non vide d'affino{\"\i}des
ouverts est un affino{\"\i}de ouvert, la collection des affino{\"\i}des ouverts
de $\pK$ satisfait la propri{\'e}t{\'e} \axiomeB{I}.  D'autre part, la
collection des affino{\"\i}des ouverts satisfait aussi la propri{\'e}t{\'e}
\axiomeB{III}, car pour chaque point $x$ dans $\pK$ et chaque $r \in
|K^*| \cap (0, 1)$, la collection
$$
\left\{ \{ z \in \pK \mid \Delta(x, z) < r \} , \{ z \in \pK \mid
  \Delta(x, z) > r/2 \} \right\}~,
$$
est un recouvrement fini de $\pK$ par des affino{\"\i}des ouverts, o{\`u} le
seul {\'e}l{\'e}ment contenant $x$ est la boule $\{ z \in \pK \mid \Delta(x,
z) < r \}$.

Le reste de ce paragraphe est d{\'e}di\'e {\`a} v\'erifier que la collection des
affino{\"\i}des ouverts v{\'e}rifie la propri{\'e}t{\'e} \axiomeB{II}.  Pour cela, on
appellera \textit{affino{\"\i}de rationnel} toute intersection finie de
boules ouvertes ou ferm{\'e}es de $\pK$.  Notons que le compl{\'e}mentaire
d'un affino{\"\i}de rationnel s'{\'e}crit de fa{\c c}on canonique comme une r{\'e}union
finie disjointe de boules ouverts ou ferm{\'e}es.  Ces boules seront
appel{\'e}es \textit{composantes du compl{\'e}mentaire} de l'affi\-no{\"\i}de.
Deux affino{\"\i}des rationnels partageant une composante de son
compl{\'e}mentaire, s'inter\-sectent.

\partn{1} Soit $C$ un recouvrement fini de $\pK$ par des affino{\"\i}des ouverts.
On note par~$\fS$ la collection finie de toutes les boules $B$ de $\pK$, tel qu'il existe un affino\"{\i}de $Y$ dans $C$ tel que $B$ ou $\pK \setminus B$ soit l'une des composantes de $\pK \setminus Y$.
Choisissons pour chaque boule ouverte $D$ dans $\fS$ un point $x_D \in D$.  Etant donn{\'e} $\eta \in
|K^*| \cap (0, 1)$, on d{\'e}signe par $A_D(\eta)$ la couronne $\{ \eta <
|z| < 1 \}$ dans une coordonn{\'e}e telle que $x_D = 0$ et $D = \{ |z| < 1
\}$.  Quitte {\`a} prendre $\eta$ plus proche de~$1$, on suppose que
$A_D(\eta)$ est contenu dans chacune des boules dans $\fS$ qu'elle
rencontre.  De plus, on suppose que si $D$ et $D'$ sont des boules
ouvertes distinctes dans $\fS$, alors les couronnes $A_D(\eta)$ et
$A_{D'}(\eta)$ sont disjointes.  Pour ce choix de $\eta$, on pose $A_D
\= A_D(\eta)$.

\partn{2} Etant donn{\'e} un point $x$ dans $\pK$, soit $Z(x)$ l'affino{\"\i}de
rationnel intersection de tous les {\'e}l{\'e}ments dans $\fS$ contenant $x$.
La collection $\fP$ de ces affino{\"\i}des rationnels est un recouvrement fini de $\pK$.
Elle est en effet une partition de $\pK$, car si $Z$ et $Z'$ sont des \'el\'ements distincts de $\fP$, alors il existe $B \in \fS$ tel que $Z \subset B$ et $Z' \subset \pK \setminus B$, d'o\`u $Z \cap Z' = \emptyset$.
Notons de plus qu'un {\'e}l{\'e}ment de~$\fP$ est contenu dans chaque {\'e}l{\'e}ment de~$C$ qu'il rencontre.

Le compl{\'e}mentaire dans $\pK$ d'un {\'e}l{\'e}ment de $\fP$ est une r{\'e}union disjointe de boules dans $\fS$.
Comme $\fP$ est une partition de $\pK$, on conclut que pour chaque boule~$D$ dans~$\fS$ il existe un unique {\'e}l{\'e}ment~$Z$ de~$\fP$ dont $D$ est l'une des composantes de~$\pK \setminus Z$.

\partn{3}
Pour chaque $Z$ dans $\fP$ on pose,
$$
Y_Z \= Z \sqcup (\bigcup_D A_D) ~,
$$
o{\`u} $D$ parcourt les boules ouvertes composantes de $\pK \setminus
Z$.  Notons que $Y_Z$ est un affino{\"\i}de ouvert.  D'autre part, notons
que tout affino{\"\i}de ouvert contenant $Z$ rencontre chacune des couronnes
$A_D$, pour $D$ comme ci-dessus.  Par cons{\'e}quent, chaque affino{\"\i}de
dans $C$ qui rencontre $Z$ contient $Z$ et chacune des
couronnes~$A_D$.  On conclut donc que chaque {\'e}l{\'e}ment de $C$ qui
rencontre $Z$ contient $Y_Z$.

\partn{4} Notons que la collection $C' \= \{ Y_Z \mid Z \in \fP \}$
est un recouvrement fini de $\pK$ par des affino{\"\i}des ouverts.  Soit
$x$ un point de $\pK$ quelconque.  On montrera que la r{\'e}union de tous
les {\'e}l{\'e}ments de $C'$ contenant $x$ est contenue dans l'un des {\'e}l{\'e}ments
de $C$.

Supposons d'abord que $x$ n'appartient {\`a} aucune des couronnes $A_D$.
Alors pour tout {\'e}l{\'e}ment $Z$ de $\fP$ tel que $Y_Z$ contient $x$, on~a
$x \in Z$.  Il existe donc un unique {\'e}l{\'e}ment $Z$ de $\fP$ avec cette
propri{\'e}t{\'e}.  Comme tout {\'e}l{\'e}ment de $C$ qui rencontre $Z$ contient $Y_Z$,
l'assertion est v{\'e}rifi{\'e}e dans ce cas.

Supposons maintenant que $x$ appartienne {\`a} l'une des couronnes $A_D$,
pour une certaine boule ouverte $D$ dans $\fS$.  Soit $Z$ l'unique
{\'e}l{\'e}ment de $\fP$ dont $D$ est l'une des composantes de son
compl{\'e}mentaire.  Alors on~a $x \in A_D \subset Y_Z$.  Comme par
d{\'e}finition les couronnes $A_D$ sont disjointes deux {\`a} deux, on conclut
que pour tout {\'e}l{\'e}ment $Z'$ de $\fP$ distinct de $Z$ et tel que
$Y_{Z'}$ contient $x$, on~a $x \in Z'$.  Par cons{\'e}quent il existe un
unique {\'e}l{\'e}ment $Z'$ de $\fP$ avec cette propri{\'e}t{\'e}.  On conclut donc
que $Y_Z$ et $Y_{Z'}$ sont les seuls {\'e}l{\'e}ments de $C'$ contenant $x$.
Soit~$Y$ un {\'e}l{\'e}ment de~$C$ qui rencontre~$Z$.
Alors~$Y$ contient $Y_Z$ et donc il rencontre~$Z'$, car $x \in Y_Z \cap Z'$.
Par cons{\'e}quent~$Y$ contient $Y_{Z'}$ et on~a $Y_Z \cup Y_{Z'} \subset Y$.

\subsection{Compactification des recouvrements finis par des affino{\"\i}des ouverts}\label{ss:berK}
On munit $\pK$ de la structure uniforme d{\'e}finie dans la
Proposition~\ref{p:structure uniforme}, lorsque $X$ est {\'e}gal {\`a} $\pK$
et $\fB$ est la collection des affino{\"\i}des ouverts.  La topologie sur
$\pK$ d\'eduite de la structure uniforme co{\"\i}ncide avec la topologie
induite par la distance $\Delta$ (Proposition~\ref{p:topologie deduite}).  On d{\'e}note par $\berK$ le s{\'e}par{\'e} compl{\'e}t{\'e} de $\pK$, voir
\S~\ref{ss:completion}.  Comme $\pK$ est s{\'e}par{\'e}, $\pK$ est
canoniquement isomorphe au sous-espace dense de $\berK$ des filtres
convergents (\S~\ref{ss:completion}).  De plus, l'espace uniforme
$\berK$ est compact (Proposition~\ref{p:le complete est compact}) et
il est donc une compactification de $\pK$.  Notons finalement que
$\berK$ est un espace uniforme r{\'e}gulier (\S~\ref{ss:completion}).

\begin{prop}\label{p:berK}
  L'espace uniforme $\berK$ est m{\'e}trisable si et seulement si le corps
  r{\'e}siduel et le groupe des valeurs de $K$ sont d{\'e}nombrables.
\end{prop}

\begin{proof}
  Supposons d'abord que le corps r{\'e}siduel et le groupe des valeurs de
  $K$ soient d{\'e}nombrables.  Alors la collection des boules ouvertes de
  $X$ et la collection des affino{\"\i}des ouvertes de $\pK$ sont
  d{\'e}nombrables.  Il s'en suit que $\berK$ est m{\'e}trisable,
  voir~\S~\ref{ss:topologie deduite}.
  
  Supposons d'autre part que le corps r{\'e}siduel de $K$ ou son groupe de
  valeurs ne soit pas d{\'e}nombrable.  Alors pour chaque $r \in |K^*|$ la
  collection $\cB_r$ des boules ouvertes de $K$ de diam{\`e}tre $r$ n'est
  pas d{\'e}nombrable.  Par cons{\'e}quent $\{ \hB \mid B \in \cB_r \}$ est
  une collection non d{\'e}nombrable de parties ouvertes de $\berK$
  (Proposition~\ref{p:base de la topologie}), qui sont deux {\`a} deux disjointes.  On
  conclut que la topologie de $\berK$ n'admet aucune base d{\'e}nombrable
  et que l'espace uniforme $\berK$ n'est pas m{\'e}trisable \cite[IX,
  \S~$4$, Th{\'e}or{\`e}me~$1$]{Bou}.
\end{proof}

\subsection{Points de $\berK$}\label{ss:points de berK}
Rappelons que les points de~$\berK$ sont les filtres de Cauchy
minimaux, voir~\S~\ref{ss:completion}.

\begin{prop}\label{p:points de berK}
  Soit $\fF$ un filtre de Cauchy minimal qui ne soit pas convergent,
  i.e. dans $\berK \setminus \pK$.  Alors la collection $\fD$ de
  toutes les boules de $K$ contenues dans $\fF$ est non vide et
  compl{\`e}tement ordonn{\'e}e par rapport {\`a} l'inclusion.
  
  De plus, pour l'intersection $B_{\fF} = \bigcap_{B \in \fD} B$ il y a
  trois cas.
\begin{enumerate}
\item[1.]  \textit{\textsf{$B_{\fF}$ est une boule ferm{\'e}e de $K$}}.  Alors la
  collection de tous les affino{\"\i}des ouverts de la forme
$$
\{ z \in K \mid r < \min \{ |z - a_j|  \mid j = 1, \ldots, k \} < r'\} ~,
$$
o{\`u} $k$ est un entier positif, $a_1, \ldots, a_k \in B_{\fF}$, $r
\in (0, \diam(B_{\fF}))$ et $r' > \diam(B_{\fF})$, est une base du
filtre $\fF$.
\item[2.]  \textit{\textsf{$B_{\fF}$ est une boule irrationnelle de $K$}}.  Alors
  la collection de toutes les couronnes de la forme $\{ z \in K \mid r
  < |z - a| < r' \}$, pour $a \in B_{\fF}$, $r \in (0,
  \diam(B_{\fF}))$ et $r' > \diam(B_{\fF})$, est une base du filtre
  $\fF$.
\item[3.]  \textit{\textsf{$B_{\fF}$ est vide}}.  Alors toute suite d{\'e}croissante
  d'{\'e}l{\'e}ments de $\fD$ dont l'intersection est vide, est une base du
  filtre $\fF$.
\end{enumerate}
Inversement, toute collection comme dans la partie~$1$ (resp.~$2$),
construite {\`a} partir d'une boule ferm{\'e}e (resp. irrationnelle) $B \=
B_\fF$ de $K$, et toute suite decroissante de boules de $K$ dont
l'intersection est vide, est une base d'un filtre de Cauchy minimal
qui n'est pas convergent.
\end{prop}

La d{\'e}monstration de cette proposition est au paragraphe
\S~\ref{ss:preuve points de ber} ci-dessous.

D'apr{\`e}s la proposition il y a~$4$ types de points de $\berK$~:
\begin{enumerate}
\item[(1)]
\textit{\textsf{Les points de $\pK$.}}
\item[(2)]
\textit{\textsf{Les points rationnels}}, d{\'e}termin{\'e}s par une boule ferm{\'e}e de~$K$.
\item[(3)]
\textit{\textsf{Les points irrationnels}}, d{\'e}termin{\'e}s par une boule irrationnelle de~$K$.
\item[(4)] \textit{\textsf{Les points singuliers}}, determin{\'e}s par une suite d{\'e}croissante de boules de~$K$, dont l'intersection est vide.
\end{enumerate}

Il est facile de voir que deux boules ferm{\'e}es ou irrationnelles
distinctes de~$K$ d{\'e}terminent des points distincts de~$\berK$.  Les
points rationnels (resp. irrationnels) de~$\berK$ sont alors en
bijection avec les boules ferm{\'e}es (resp. irrationnelles) de~$K$.

\subsection{D{\'e}monstration de la Proposition~\ref{p:points de berK}}\label{ss:preuve points de ber}
La d{\'e}monstration de la Proposition~\ref{p:points de berK} est
ci-dessous~; elle s'appuit sur les lemmes~\ref{l:filtre de Cauchy
  minimal} et~\ref{l:description de berK}.
\begin{lemm}\label{l:filtre de Cauchy minimal}
  Pour qu'un filtre de Cauchy $\fF$ soit minimal, il suffit qu'il
  admette une base form{\'e}e d'affino{\"\i}des ouverts et que tout {\'e}l{\'e}ment de
  $\fF$ contienne un affino{\"\i}de ferm{\'e} appartenant {\`a} $\fF$.
\end{lemm}
\begin{proof}
  Soit $\fF$ un filtre de Cauchy satisfaisant la propri{\'e}t{\'e} d{\'e}crite et
  soit $\fF'$ le filtre de Cauchy minimal contenu dans $\fF$.  Comme
  par hypoth{\`e}se $\fF$ admet une base form{\'e}e d'affino{\"\i}des ouverts, il
  suffit de montrer que pour tout affino{\"\i}de ouvert $Y$ contenu dans
  $\fF$, on~a $Y \in \fF'$.
  
  Par hypoth{\`e}se il existe un affino{\"\i}de ferm{\'e} $X \in \fF$ contenu dans
  $Y$.  Soient $B_1, \ldots, B_k$ les composantes de $\pK \setminus
  X$.  Alors $C = \{ Y, B_1, \ldots, B_k \}$ est un recouvrement fini
  de $\pK$ par affino{\"\i}des ouverts.  Le Lemme~\ref{l:filtres de Cauchy}
  implique qu'au moins l'un des {\'e}l{\'e}ments de $C$ appartient {\`a} $\fF$.
  Mais, comme $X \in \fF$, aucune des boules $B_i$ n'appartient {\`a}
  $\fF$.  Comme $\fF' \subset \fF$ aucune des boules $B_i$
  n'appartient {\`a} $\fF$ et alors le Lemme~\ref{l:filtres de Cauchy}
  implique qu'on a $Y \in \fF'$.  Ceci termine la preuve du lemme.
\end{proof}
\begin{lemm}\label{l:description de berK}
  Soit $\fB$ une collection d'affino{\"\i}des ouverts, d{\'e}crite {\`a} partir
  d'une boule ferm{\'e}e (resp. irrationnelle) $B$ de $K$, comme dans la
  partie~$1$ (resp.~$2$) de la proposition, ou soit $\fB$ une suite
  d{\'e}croissante de boules de $K$ dont l'intersection est vide.  Alors
  $\fB$ est une base d'un  filtre de Cauchy minimal qui n'est pas
  convergent.
\end{lemm}
\begin{proof}
  Il est facile de voir que toute intersection finie d'{\'e}l{\'e}ments de
  $\fB$ contient un {\'e}l{\'e}ment de $\fB$.  Par cons{\'e}quent $\fB$ engendre
  un filtre, qu'on note par $\fF$.  Comme l'intersection des {\'e}l{\'e}ments
  de $\fB$ est vide, le filtre $\fF$ n'est pas convergent.  Il est
  facile de voir que $\fB$, et donc $\fF$, v{\'e}rifie la propri{\'e}t{\'e}
  d{\'e}crite dans le Lemme~\ref{l:filtre de Cauchy minimal}.  Il suffit
  alors de montrer que le filtre $\fF$ est de Cauchy.  D'apr{\`e}s le
  Lemme~\ref{l:filtres de Cauchy} il faut montrer que $\fF$ contient
  un {\'e}l{\'e}ment de chaque recouvrement dans $\fC$.

Soit alors $C \in \fC$.
Il y a trois cas.

\case{1}{$\fB$ est d{\'e}finie {\`a} partir d'une boule ferm{\'e}e $B$.}
Apr{\`e}s changement de coordonn{\'e}e affine, on suppose $B = \OK = \{ |z|
\le 1 \}$.  Comme par hypoth{\`e}se le corps r{\'e}siduel de $K$ est infini,
il existe un {\'e}l{\'e}ment $Y$ de $C$ qui rencontre au moins~$2$ classes
r{\'e}siduelles.  Chaque composante de $\pK \setminus Y$ est alors, soit
disjointe de $\OK$, soit contenue dans $\OK$.  En particulier il
existe $r' > 1$ tel que $\{ 1 < |z| < r' \} \subset Y$.  On d{\'e}signe
par $B_1, \ldots, B_k$ les composantes de $\pK \setminus Y$ contenues
dans $\OK$.  Comme $Y$ rencontre au moins~$2$ classes r{\'e}siduelles,
on~a $\max \{ \diam(B_j) \} < 1$.  On choisit $r \in (\max\{
\diam(B_j) \}, 1)$ et pour chaque $j = 1, \ldots, k$ on choisit un
point $a_j \in B_j$.  Alors l'affino{\"\i}de ouvert $\{ z \in K \mid r <
\min \{ |z - a_j| \} < r'\}$ est un {\'e}l{\'e}ment de $\fB$ contenu dans $Y$.
Il s'en suit que $Y \in \fF$.

\case{2}{$\fB$ est d{\'e}finie {\`a} partir d'une boule irrationnelle $B$.}  Apr{\`e}s changement de coordonn{\'e}e affine, on suppose que $B$ est
de la forme $\{ |z| < r_0 \}$, avec $r_0 \not \in |K^*|$.  Notons que
tout affino{\"\i}de ouvert contenu dans $B$ (resp. $\pK \setminus B$) est
contenu dans une boule ouverte contenue dans $B$ (resp. $\pK \setminus
B$).  Par cons{\'e}quent il existe $Y \in C$ qui rencontre $B$ et $\pK
\setminus B$.  Alors il est facile de voir que toute composante de
$\pK \setminus Y$ est contenue dans $B$ ou dans $\pK \setminus B$.  Il
existe donc $r \in (0, r_0)$ et $r' > r_0$ dans $|K^*|$, tels que la
couronne $\{ z \in K \mid r < |z| < r' \} \in \fB$ soit contenue dans
$Y$.  On a alors $Y \in \fF$.

\case{3}{$\fB$ est une suite d{\'e}croissante de boules de $K$ dont l'intersection est vide.}
Consid{\'e}rons une boule $D$ de $\pK$
quelconque.  Si $D \subset K$, alors on choisit un point $x$ dans $D$
et sinon, alors on choisit un point $x$ dans $K \setminus D$.  Comme
$\bigcap_{B \in \fB} B = \emptyset$ il existe une boule $B_0 \in \fB$
contenue dans $K \setminus \{ x \}$.  Il est facile de voir que dans
ce cas on~a $B_0 \subset D$ ou $B_0 \subset \pK \setminus D$.

Consid{\'e}rons maintenant un affino{\"\i}de ouvert $Y = \bigcap_{j = 1, \ldots,
  j} D_j$ quelconque.  Si l'on applique le raisonnement pr{\'e}c{\'e}dent {\`a}
chacune des boules $D_i$, on conclut qu'il existe une boule $B_0$ dans
$\fB$ telle que $B_0 \subset Y$ ou telle que $B_0 \subset \pK
\setminus Y$.  Il s'en suit qu'il existe $Y \in C$ et $B \in \fB$ tel
que $B \subset Y$.  On a alors $Y \in \fF$.
\end{proof}

\begin{proof}[D{\'e}monstration de la Proposition~\ref{p:points de berK}.]
  La derni{\`e}re assertion de la proposition est donn{\'e}e par le
  Lemme~\ref{l:description de berK}.
  
  L'ensemble $\fD$ est compl{\`e}tement ordonn{\'e} par rapport {\`a} l'inclusion,
  car deux boules quelconques de $K$ qui s'intersectent sont reli\'ees
  par l'inclusion.  Montrons maintenant que $\fD$ n'est pas vide.
  Comme par hypoth{\`e}se $\fF$ n'est pas convergent, on~a en particulier
  $\fF \neq \infty$ et par cons{\'e}quent il existe $r > 0$ tel que $\{
  |z| > r \} \cup \{ \infty \} \not\in \fF$.  Le Lemme~\ref{l:filtres
    de Cauchy} entra{\^\i}ne donc que $\{ |z| < 2 \cdot r \} \in \fF$.
  
  Supposons d'abord que l'ensemble $B_\fF$ est vide et soit $\fB$ une
  suite d{\'e}croissante de boules dans $\fD$ dont l'intersection est
  vide.  Le Lemme~\ref{l:description de berK} implique alors que $\fB$
  est une base d'un filtre de Cauchy minimal.  Ce filtre est contenu
  dans $\fF$ et par minimalit{\'e} il est {\'e}gal {\`a}~$\fF$.  Il s'en suit que
  $\fB$ est une base de $\fF$.
  
  Supposons maintenant que l'ensemble $B_\fF$ est non vide.  Alors
  pour chaque $a_0 \in B_\fF$ on~a
$$
B_\fF = \{ z \in K \mid |z - a_0| \le r_0 \},
$$
pour un certain $r_0 \ge 0$.  Comme par hypoth{\`e}se $\fF$ n'est pas
convergent, on~a $r_0 > 0$ et donc $B_\fF$ est une boule ferm{\'e}e ou
irrationnelle de~$K$.  Par le Lemme~\ref{l:description de berK}, la
collection d'affino{\"\i}des ouverts $\fB_{\fF}$, d{\'e}crite dans la
partie~$1$ lorsque $r_0 \in |K^*|$ et dans la partie~$2$ lorsque $r_0
\not \in |K^*|$, est une base d'un filtre de Cauchy minimal.  Pour
montrer que $\fB_{\fF}$ est une base de $\fF$, il suffit alors de
montrer que chaque {\'e}l{\'e}ment de $\fB_{\fF}$ appartient {\`a} $\fF$.  Il
suffit donc de montrer que pour chaque $a \in B_\fF$, chaque $r \in
(0, r_0)$ et chaque $r' > r_0$, la couronne $\{ r < |z - a| < r' \}$
appartient {\`a} $\fF$.  Par d{\'e}finition de $r_0$, la boule $\{ |z - a| <
r' \}$ appartient {\`a} $\fF$ et pour chaque $\rho \in (0, r_0)$ la boule
$\{ |z - a| < \rho \}$ n'appartient pas {\`a} $\fF$.  Si l'on choisit
$\rho \in (r, r_0)$ alors,
$$
\left \{ \{ |z - a| > r \} \cup \{ \infty \}, \{ |z - a| < \rho \} \right\}~,
$$
est un recouvrement de $\pK$ dans $\fC$.  Le Lemme~\ref{l:filtres
  de Cauchy} implique alors qu'on a $\{ |z - a| > r \} \cup \{ \infty
\} \in \fF$ et donc $\{ r < |z - a| < r' \} \in \fF$.  Ceci termine la
preuve de la proposition.
\end{proof}
\subsection{Carat{\'e}risation des filtres de Cauchy minimaux}\label{ss:caracteristique minimalite}
\

\begin{coro}\label{c:caracteristique minimalite}
  Un filtre de Cauchy $\fF$ est minimal si et seulement si il admet
  une base form{\'e}e d'affino{\"\i}des ouverts et si tout {\'e}l{\'e}ment de $\fF$
  contient un affino{\"\i}de ferm{\'e} appartennant {\`a} $\fF$.
\end{coro}
\begin{proof}
  Il est clair que tout filtre convergent satisfait la propri{\'e}t{\'e}
  d{\'e}crite.  D'autre part, on voit facilement de la
  Proposition~\ref{p:points de berK} que tout filtre de Cauchy minimal
  qui n'est pas convergent satisfait la propri{\'e}t{\'e} d{\'e}crite.
  L'implication inverse est donn{\'e}e par le Lemme~\ref{l:filtre de
    Cauchy minimal}.
\end{proof} 
\subsection{Lien avec la droite projective de Berkovich}\label{ss:lien}
Par d\'efinition, l'espace analytique de $K$ au sens de Berkvoich est \'egal \`a l'espace de toutes les seminormes multiplicatives et continues dans $K[T]$, muni de la plus petite topologie que rend toutes ces seminormes continues~\cite{Berrouge}.
On obtient l'espace analytique de $\pK$ par recollement de deux copies de l'espace analytique de~$K$, de la fa\c{c}on usuelle.

Le lien entre cet espace analytique avec l'espace uniforme $\berK$ qu'on a d\'ecrit ici, et alors tr\`es explicite.
Chaque filtre de Cauchy minimal $\fF$ distinct de~$\infty$ d{\'e}finit une semi-norme multiplicative et continue $| \cdot |_\fF$ dans $K[T]$, par
$$
|P|_\fF \= \lim_{\fF} |P| ~.
$$
Inversement, toute semi-norme multiplicative de $K[T]$, qui n'est
pas de la forme $P \mapsto |P(x)|$ pour un certain point $x \in K$,
est de cette forme.
De plus, l'application $\fF \mapsto | \cdot |_\fF$ est un homeomorphisme entre $\berK \setminus \{ \infty \}$ et l'espace analytique de~$K$.

Notons d'autre part que la description des quatre types de points de $\berK$ en~\S~\ref{ss:points de berK}, est analogue \`a la description dans~\cite[p.5]{Berrouge}.
Les points de $\pK$ sont appelles points de \og{}type~\textrm{(1)}\fg{} dans~\cite{Berrouge}, les points qu'on appelle ici \og{}points rationels\fg{} correspondent aux points de \og{}type~\textrm{(2)}\fg{}, les points \og{}irrationnels\fg{} aux points de \og{}type~\textrm{(3)}\fg{}, et les points \og{}singuliers\fg{} aux points de \og{}type~\textrm{(4)}\fg{}.

\newpage
\section{G{\'e}om{\'e}trie de $\berK$}\label{s:geometrie de berK}
Comme dans la section pr{\'e}c{\'e}dente, on fixe un corps $K$ muni d'une
norme ultram{\'e}trique $| \cdot |$ pour laquelle $K$ est complet et on
suppose que le corps r{\'e}siduel de $K$ est infini et que le groupe des
valeurs de~$K$ n'est pas discret.
\subsection{Boules et affino\"{\i}des de $\berK$}\label{ss:boules et affinoides de berK}
Etant donn{\'e}e une boule ouverte ou irrationnelle (resp. ferm{\'e}e ou
irrationnelle) $B$ de $\pK$, on appelle $\hB$ (resp. $\ov{B}$)
\textit{boule ouverte} (resp. \textit{ferm{\'e}e}) \textit{de $\berK$}.
Notons que le compl{\'e}mentaire d'une boule ouverte (resp. ferm{\'e}e) de
$\berK$ est une boule ferm{\'e}e (resp. ouverte) de $\berK$.  De plus,
toute boule ouverte (resp. ferm{\'e}e) de $\berK$ est un ensemble ouvert
(resp. ferm{\'e}) de $\berK$ (Proposition~\ref{p:base de la topologie}).

Etant donn{\'e} un affino\"{\i}de ouvert (resp. ferm{\'e}) $Y$ de $\pK$, on
appelle $\hY$ (resp. $\ov{Y}$) \textit{affino\"{\i}de ouvert} (resp.
\textit{ferm{\'e}}) \textit{de $\berK$}.  Le lemme ci-dessous implique que
toute r{\'e}union finie avec intersection non vide et toute intersection
non vide d'affino\"{\i}des ouverts (resp. ferm{\'e}s) de $\berK$ est un
affino\"{\i}de ouvert (resp. ferm{\'e}) de $\berK$.  En particulier tout
affino\"{\i}de ouvert (resp. ferm{\'e}) de $\berK$ s'{\'e}crit comme une
intersection finie non vide de boules ouvertes (resp. ferm{\'e}es) de
$\berK$.

Notons que tout affino\"{\i}de ouvert (resp. ferm{\'e}) est ouvert (resp. ferm{\'e})
dans $\berK$ (Proposition~\ref{p:base de la topologie}).  De plus, la collection des
affino{\"\i}des ouverts de $\berK$ est une base de la topologie de $\berK$
(Proposition~\ref{p:base de la topologie}).
\begin{lemm}\label{l:boules et affinoides de berK}
Soient $Y$ et $Y'$ des affino\"{\i}des ouverts (resp. ferm{\'e}s) de $\pK$.
Alors on~a
$$
\widehat{Y \cap Y'} = \hY \cap \widehat{Y'}
\ \text{ et } \
\widehat{Y \cup Y'} = \hY \cup \widehat{Y'}
$$
$$
\text{(resp. $
\ov{Y \cap Y'} = \ov{Y} \cap \ov{Y'}
\ \text{ et } \
\ov{Y \cup Y'} = \ov{Y} \cup \ov{Y'}$)}~.
$$
\end{lemm}
\begin{proof}
  A l'aide de la Proposition~\ref{p:base de la topologie}, les assertions concernant les
  affi\-no{\"\i}des ferm{\'e}s se d{\'e}duisent de celles concernant les affino{\"\i}des
  ouverts.  Par~\eqref{e-core} il suffit donc de montrer que si $Y$ et
  $Y'$ sont des affino{\"\i}des ouverts, alors on~a $\widehat{Y \cup Y'}
  \subset \hY \cup \widehat{Y'}$.

Soit $\fF \in \widehat{Y \cup Y'}$.  Il existe alors un affino{\"\i}de
ferm{\'e} $X \in \fF$ contenu dans $Y \cup Y'$
(Corollaire~\ref{c:caracteristique minimalite}).  Par cons{\'e}quent, la
collection $C$ form{\'e}e de $Y$, $Y'$ et des boules ouvertes composantes
de $\pK \setminus X$, est un recouvrement fini de $\pK$.  Le
Lemme~\ref{l:filtres de Cauchy} implique donc que $\fF$ contient l'un
des {\'e}l{\'e}ments de $C$.  Comme $X \in \fF$, le filtre $\fF$ ne contient
aucune composante de $\pK \setminus X$~; $\fF$ contient alors $Y$ ou
$Y'$.
\end{proof}

\subsection{Le point de $\berK$ determin{\'e} par une boule de $\pK$}\label{ss:point d'une boule}
Rappelons que pour un point rationnel ou irrationnel $\fF$ de $\berK$,
on d{\'e}signe par $B_{\fF}$ la boule ferm{\'e}e ou irrationnelle de $K$
d\'etermin{\'e}e par $\fF$, voir~\S~\ref{ss:points de berK}.

Soit $B$ une boule de $\pK$ et soit $B' = \pK \setminus B$ la boule
compl{\'e}mentaire.  Notons qu'on a, soit $B \subset K$, soit $B' \subset
K$.  Lorsque $B \subset K$ (resp. $B' \subset K$), soit $\fF$ le point
de $\berK$ determin{\'e} par l'unique boule ferm{\'e}e ou irrationnelle de $K$
contenant $B$ (resp. $B'$) et du m{\^e}me diam{\`e}tre que $B$ (resp. $B'$).
Notons que si $B$ est ouverte ou ferm{\'e}e (resp. irrationnelle), alors
$\fF$ est rationnel (resp. irrationnel).

De cette fa{\c c}on chaque boule $B$ de $\pK$ d\'etermine un point de
$\berK$.  Notons que deux boules compl\'ementaires d\'eterminent le m{\^e}me
point.  Le lemme suivant est une cons{\'e}quence imm{\'e}diate de la
Proposition~\ref{p:points de berK}.

\begin{lemm}\label{l:point d'une boule}
  Soit $B$ une boule de $\pK$ et soit $\fF$ le point de $\berK$
  d\'etermin{\'e} par~$B$.  Alors $B \not \in \fF$ et $B$ rencontre tout
  {\'e}l{\'e}ment de $\fF$.
\end{lemm}
\subsection{Boules de $\pK$ associ{\'e}es {\`a} un point de $\berK$}\label{ss:boules d'un point}
Soit $\fF$ un point rationnel ou irrationnel de $\berK$.  On dira
qu'une boule $B$ de $\pK$ est \textit{associ{\'e}e {\`a}} $\fF$, si $\fF$ est
le point de $\berK$ d\'etermin{\'e} par $B$.  On d{\'e}signe par $\pi(\fF)$ la
collection de toutes les boules ouvertes ou irrationnelles de $\pK$
associ{\'e}es {\`a} $\fF$.  Notons que lorsque $\fF$ est rationnel (resp.
irrationnel) tous les {\'e}l{\'e}ments de $\pi(\fF)$ sont des boules ouvertes
(resp. irrationnelles) de $\pK$.

Il facile de voir que la boule $B(\infty) = \pK \setminus B_{\fF}$
appartient {\`a} $\pi(\fF)$.  De plus, toute boule appartenant {\`a}
$\pi(\fF)$ et distincte de $B(\infty)$ est contenue dans $B_{\fF}$ et
tout point de $B_{\fF}$ est contenu dans l'un des {\'e}l{\'e}ments de
$\pi(\fF)$.  La propri{\'e}t{\'e} ultram{\'e}trique implique alors que les
{\'e}l{\'e}ments de $\pi(\fF)$ distincts de $B(\infty)$ sont deux {\`a} deux
disjoints.  On conclut que les {\'e}l{\'e}ments de $\pi(\fF)$ forment une
partition de $\pK$~; en particulier on~a
$$
\pK = \sqcup_{B \in \pi(\fF)} B ~.
$$

Lorsque $\fF$ est irrationnel on~a $\pi(\fF) = \{ B_{\fF}, \pK
\setminus B_{\fF} \}$.

Dans le cas o{\`u} $B_{\fF}$ est la boule unit{\'e} $\cO_K$, les {\'e}l{\'e}ments de
$\pi(\fF)$ co{\"\i}ncident avec les classes r{\'e}siduelles~; i.e. avec les
fibres de la projection de $\pK$ vers $\ptK$.  On a donc un
param{\'e}trage de $\pi(\fF)$ par $\ptK$.  Lorsque $\fF$ est un point
rationnel quelconque, on peut se ram{\`e}ner au cas pr{\'e}c{\'e}dent par un
changement affine de coordonn{\'e}es.  Alors on~a un param{\'e}trage de
$\pi(\fF)$ par $\ptK$, qui est unique modulo une transformation affine
de $\ptK$.
\subsection{La partition de $\berK$ d\'etermin{\'e}e par un point}\label{ss:partition d'un point}
\ 

\begin{prop}\label{p:partition d'un point}
  Soit $\fF$ un point rationnel ou irrationnel de $\berK$.  Alors les
  ensembles de la forme $\hB$, avec $B \in \pi(\fF)$, sont disjoints
  deux {\`a} deux et on~a
$$
\berK \setminus \{ \fF \} = \sqcup_{B \in \pi(\fF)} \hB~.
$$
\end{prop}
\begin{proof}
  Comme les {\'e}l{\'e}ments de $\pi(\fF)$ sont disjoints deux {\`a} deux, les
  ensembles de la forme $\hB$, avec $B \in \pi(\fF)$, sont disjoints
  deux {\`a} deux.  Par le Lemme~\ref{l:point d'une boule} pour chaque $B
  \in \pi(\fF)$ on~a $B \not \in \fF$~; c'est-{\`a}-dire $\fF \not \in
  \hB$.  Il reste \`a montrer que pour tout point $\fF'$ de $\berK$
  distinct de $\fF$ il existe $B \in \pi (\fF)$ tel que $B \in \fF'$.
  Soit $Y \in \fF$ un affino{\"\i}de ouvert tel que $Y \not \in \fF'$.  Par
  la Proposition~\ref{p:points de berK} il est facile de voir qu'il
  n'y a qu'un nombre fini d'{\'e}l{\'e}ments $B_1, \ldots, B_n$ de $\pi(\fF)$
  qui ne sont pas contenus dans $Y$.  Alors $\{ Y, B_1, \ldots, B_n
  \}$ est un recouvrement fini de $\pK$ par d'affino{\"\i}des ouverts.  Par
  cons{\'e}quent $\fF'$ contient l'un des {\'e}l{\'e}ments de ce recouvrement.
  Comme $Y \not \in \fF'$ on conclut qu'il existe $j = 1, \ldots n$
  tel que $B_j \in \fF'$.  On a donc $\fF' \in \widehat{B_j}$.
\end{proof}
\subsection{Fronti{\`e}re d'une boule}\label{ss:frontiere d'une boule}
Rappelons que pour une partie $Y$ de $\pK$ la fermeture topologique
$\ov{\hY}$ de $\hY$ est {\'e}gale {\`a} $\ov{Y}$, voir~\S~\ref{ss:topologie du complete}.

\begin{lemm}\label{l:frontiere d'une boule}
  Soit $B$ une boule de $\pK$ et soit $\fF$ le point de $\berK$
  d\'etermin{\'e}e par $B$.  Alors on~a
$$
\fF \not \in \hB
\ \text{ et } \
\ov{B} = \hB \sqcup \{ \fF \} ~.
$$
En particulier l'int\'erieur de $\ov{B}$ est {\'e}gal {\`a} $\hB$ et les
ensembles $\ov{B}$ et $\hB$ ont la m{\^e}me fronti{\`e}re topologique, {\'e}gale {\`a}
$\{ \fF \}$.
\end{lemm}
\begin{proof}
  Le Lemme~\ref{l:point d'une boule} implique qu'on a $\hB \cap \{ \fF
  \} = \emptyset$ et $\fF \in \ov{B}$.

Lorsque la boule $B$ est ouverte ou irrationnelle, la
Proposition~\ref{p:partition d'un point} implique que le
compl{\'e}mentaire de $\hB \cup \{ \fF \}$ dans $\berK$ est {\'e}gal {\`a}
l'ensemble ouvert $\bigcup_{B' \in \pi(\fF), \; B' \neq B} \widehat{B'}$.
Dans ce cas on~a donc $\ov{B} = \hB \sqcup \{ \fF \}$.

Lorsque la boule $B$ est ferm{\'e}e, l'ensemble $B' = \pK \setminus B$ est
une boule ouverte et le point de $\berK$ d\'etermin{\'e} par $B'$ est $\fF$.
Par ce qui pr{\'e}c\`ede on~a donc $\ov{B'} = \widehat{B'} \sqcup \{ \fF \}$
et la Proposition~\ref{p:base de la topologie} implique qu'on a $\hB = \berK \setminus
\ov{B'} = \ov{B} \setminus \{ \fF \}$.
\end{proof}
\subsection{Fronti{\`e}re d'un affino{\"\i}de}\label{frontiere d'un affinoide}
\

\begin{lemm}\label{l:frontiere d'un affinoide}
  Soit $Y$ un affino{\"\i}de ouvert ou ferm{\'e} de $\pK$ et soient $B_1,
  \ldots, B_n$ les boules de $\pK$ composantes de $\pK \setminus Y$.
  Pour $j = 1, \ldots, n$ on d{\'e}signe par $\fF_j$ le point de $\berK$
  d\'etermin{\'e} par la boule $B_j$.  Alors on~a
$$
\hY \cap \{ \fF_1, \ldots, \fF_n \} = \emptyset
\ \text{ et } \
\ov{Y} = \hY \sqcup \{ \fF_1, \ldots, \fF_n \} ~.
$$
En particulier l'int\'erieur de $\ov{Y}$ est {\'e}gal {\`a} $\hY$ et les
ensembles $\hY$ et $\ov{Y}$ ont la m{\^e}me fronti{\`e}re topologique, {\'e}gale {\`a}
$\{ \fF_1, \ldots, \fF_n \}$.
\end{lemm}
\begin{proof}
  Pour chaque $j$ on pose $D_j \= \pK \setminus B_j$.  Alors on~a $Y =
  D_1 \cap \ldots \cap D_n$ et le Lemme~\ref{l:point d'une boule}
  implique que $\hY$ est disjoint de $\{ \fF_1, \ldots, \fF_n \}$.
  
  Lorsque l'affino{\"\i}de $Y$ est ouvert, les boules $D_j$ sont ouvertes
  et les lemmes~\ref{l:boules et affinoides de berK}
  et~\ref{l:frontiere d'une boule} impliquent
$$
\ov{Y}
=
\ov{D_1} \cap \ldots \cap \ov{D_n}
=
\hY \sqcup \{ \fF_1, \ldots, \fF_n \}~.
$$

Supposons maintenant que l'affino{\"\i}de $Y$ est ferm{\'e}.  Alors chacune des
boules $B_j$ est ouverte et les lemmes~\ref{l:boules et affinoides de
  berK} et~\ref{l:frontiere d'une boule} impliquent qu'on a
$$
\ov{\pK \setminus Y} = \ov{B_1} \cup \ldots \cup \ov{B_n} = (\hB_1
\cup \ldots \cup \hB_n) \cup \{ \fF_1, \ldots, \fF_n \}~.
$$
La Proposition~\ref{p:base de la topologie} implique alors qu'on a
$$
\hY = (\berK \setminus (\hB_1 \cup \ldots \cup \hB_n)) \setminus \{
\fF_1, \ldots, \fF_n \} = \ov{Y} \setminus \{ \fF_1, \ldots, \fF_n
\}~.
$$
Comme les ensembles $\hY$ et $\{ \fF_1, \ldots, \fF_n \}$ sont
disjoints, on obtient l'assertion d{\'e}sir{\'e}e.
\end{proof}

\newpage
\section{L'arbre associ{\'e} {\`a} un espace ultram{\'e}trique}\label{s:l'arbre d'un espace}
Dans cette section on~associe {\`a} chaque espace ultram{\'e}trique $(X,
\distX)$ un arbr{\'e} r{\'e}el $\arbX$ o{\`u} $X$ se plonge isom{\'e}triquement.  On
pourra consulter~\cite{Hug} ou~\cite{Le} pour les propri{\'e}t{\'e}s
fonctorielles de quelques variantes de cette construction.  Dans cette
section on suit la notation de~\cite[\S\S 3, 4]{elements}.
Par simplicit\'e on se resteindra au cas des espaces m\'etriques complets.
\subsection{Espace des boules}\label{ss:espace des boules}
Fixons un espace ultram{\'e}trique complet $(X, \distX)$ de diam{\`e}tre $\diam(X) \in
[0, + \infty]$ et posons
$$
I =
\begin{cases}
[0, \diam(X)] & \text{lorsque $\diam(X) < + \infty$}, \\
[0, + \infty) & \text{lorsque $\diam(X) = + \infty$}.
\end{cases}
$$
Consid{\'e}rons la relation d'{\'e}quivalence $\sim$ sur $X \times I$ d{\'e}finie par
$$
(x, r) \sim (x', r')
\ \text{ si et seulement si } \
r = r' \text{ et } r \ge \dist(x, x')~,
$$
voir figure.
On d{\'e}signe par $\arbXR$ le quotient $X \times I /
\sim$ et pour chaque {\'e}l{\'e}ment $(x, r)$ de $X \times I$ on d{\'e}signe par
$[x, r]$ le point de $\arbXR$ repr{\'e}sent{\'e} par $(x, r)$.
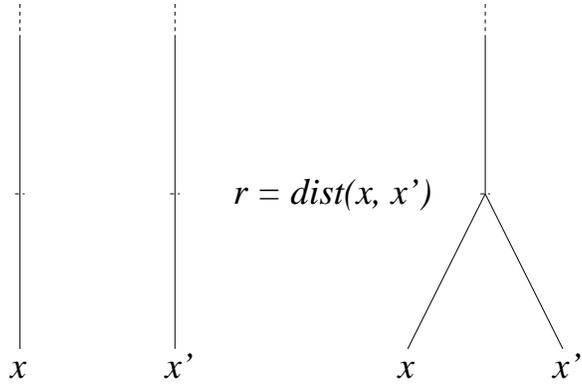
\begin{figure}[hbt]
\begin{center}
\includegraphics[width = 3in]{R.eps}
\end{center}
\end{figure}

Notons que l'application $x \mapsto [x, 0]$ de $X$ {\`a} $\arbXR$ est une
bijection sur son image. On identifiera $X$ {\`a} son image dans $\arbXR$.
D'autre part, notons que lorsque $X$ est de diam{\`e}tre fini, l'ensemble
$X \times \{ \diam(X) \}$ est une classe d'{\'e}quivalence de~$\sim$.

Notons que pour chaque $(x, r) \in X \times I$ l'ensemble
$$
\{ w \in X \mid \distX(w, x) \le r \}
$$
ne d{\'e}pend que de la classe d'{\'e}quivalence de $(x, r)$, i.e. de $[x,
r] \in \arbKR$.  On notera cet ensemble par $B_{[x, r]}$.  Notons que
lorsque $r > 0$ cet ensemble est une boule de $X$.

Pour $[x, r] \in \arbKR$ on pose
$$
\diam( [x, r] ) \= r ~.
$$
La fonction $\diam : \arbXR \to I$ ainsi d{\'e}finie s'annule
pr{\'e}cis{\'e}ment sur $X$.
\subsection{Ordre partiel $\preccurlyeq$}\label{ss:ordre partiel}
On d{\'e}signe par  $\preccurlyeq$  l'ordre partiel sur $\arbXR$ d{\'e}finie par,
$$
\cS \preccurlyeq \cS'
\ \text{ si et seulment si } \
B_{\cS} \subset B_{\cS'} ~.
$$
Si l'on pose $\cS = [x, r]$ et $\cS' = [x', r']$, alors les points
$\cS$ et $\cS'$ sont reli{\'e}s par $\preccurlyeq$ si et seulement si
$\distX(x, x') \le \max \{ r, r' \}$ et dans ce cas on~a $\cS
\preccurlyeq \cS'$ si et seulement si $r \le r'$.  En particulier $\cS
\preccurlyeq \cS'$ implique $\diam(\cS) \le \diam(\cS')$.
\subsection{Fonctions  $\cdot \vee \cdot$  et  $\sup \{ \cdot, \cdot \}$}
Notons que pour $\cS = [x, r]$ et $\cS' = [x', r']$ dans $\arbXR$, on~a
$$
[x, \max \{r, r', \distX(x, x') \}] = [x', \max \{ r, r', \distX(x,
x') \}] ~.
$$
On notera ce point de $\arbXR$ par $\cS \vee \cS'$.
Il est caract\'eris{\'e} comme le plus petit point qui est {\`a} la fois plus grand que $\cS$ et que $\cS'$, par rapport {\`a} l'ordre partiel  $\preccurlyeq$~.
D'autre part, on~a
\begin{equation}\label{e-definition ordre}
\cS \preccurlyeq \cS'
\ \text{ si et seulement si } \
\cS \vee \cS' = \cS' ~.
\end{equation}
En particulier, pour chaque point $\cS$ de $\arbXR$ on~a $\cS \vee
\cS = \cS$.

Pour $\cS$ et $\cS'$ comme avant on pose,
\begin{eqnarray*}
\sup\{ \cS, \cS' \}
& \= &
\diam (\cS \vee \cS') \\
& = &
\max \{ r, r', \distX(x, x') \} ~.
\end{eqnarray*}
Notons qu'on a
$$
\max \{ \diam(\cS), \diam(\cS') \}
\le \sup \{ \cS, \cS' \} ~,
$$
et pour $\cS_0, \cS_1 \in \arbXR$,
\begin{equation}\label{ultrametrique}
\sup\{ \cS_0, \cS_1 \}
\le \max\{ \sup\{ \cS_0, \cS \}, \sup\{ \cS, \cS_1\} \} ~.
\end{equation}
Notons finalement qu'on a $\sup \{ \cS, \cS \} = \diam(\cS)$.
\subsection{Segments}\label{ss:segments}
Etant donn{\'e}s $\cS = [x, r]$, $\cS' = [x', r'] \in \arbXR$ on pose,
\begin{eqnarray*}
[ \cS, \cS' ]
& \= &
\{ \tcS \in \arbXR \mid \cS \preccurlyeq \tcS \preccurlyeq \cS \vee \cS'
\ \text{ ou } \
\cS' \preccurlyeq \tcS \preccurlyeq \cS \vee \cS' \}
\\ & = &
\{ [x, r''] \mid r'' \in [r, \max \{ r, r', \distX(x, x') \}] \}
\\ & & \cup
\{ [x', r''] \mid r'' \in [r, \max \{ r, r', \distX(x, x') \} ] \} ~.
\end{eqnarray*}
Notons qu'on a $[\cS', \cS] = [\cS, \cS']$.
On utilisera la notation usuelle d'intervalles~:
$$
[\cS, \cS') = (\cS', \cS] = [\cS, \cS'] \setminus \{ \cS' \}
\text{ et }
(\cS, \cS') = (\cS', \cS) = [\cS, \cS') \setminus \{ \cS \}.
$$
De plus on dira qu'un point $\tcS$ de $\arbXR$ est \textit{entre} deux
points distincts $\cS$ et $\cS'$ de $\arbXR$ si $\tcS \in (\cS,
\cS')$, et on dira qu'une partie de $\arbXR$ est \textit{convexe} si et
seulement si pour toute paire de points $\cS$ et $\cS'$ appartenant {\`a}
cet ensemble, le segment $[\cS, \cS']$ appartient {\`a} cet ensemble.

Soient $\cS = [x, r]$ et $\cS' = [x', r']$ dans $\arbXR$.  Alors il
est facile de voir que l'ensemble $[\cS, \cS']$ est l'union disjointe
de $\cS \vee \cS'$ et des ensembles
$$
[\cS, \cS \vee \cS')
=
\{  [x, r''] \mid r'' \in [r, \max \{ r, r', \distX(x, x') \})  \} ~,
$$
$$
[\cS', \cS \vee \cS')
=
\{  [x', r''] \mid r'' \in [r', \max \{ r, r', \distX(x, x') \})  \}~.
$$
D'autre part, notons que la fonction  $\diam$  induit une bijection croi\-ssante entre $[\cS, \cS \vee \cS']$ (resp. $[\cS', \cS \vee \cS']$) et l'intervalle $[\diam(\cS), \sup \{ \cS, \cS' \} ]$ (resp. $[\diam(\cS'), \sup \{ \cS, \cS' \}]$) de $\R$.

\subsection{Distance $\distAX$}\label{distance delta}
Etant donn{\'e} deux points $\cS = [x, r]$ et $\cS' = [x', r']$ de $\arbXR$, on pose
\begin{eqnarray*}
\distAX(\cS, \cS')
& \= & \sup \{ \cS, \cS' \} - \tfrac{1}{2} (\diam(\cS) + \diam(\cS'))
\\ & = &
\max \{ r, r', \distX(x, x') \} - \tfrac{1}{2} (r + r') ~.
\end{eqnarray*}
On v{\'e}rifie ais\'ement que $\distAX$ d{\'e}finit une distance sur $\arbXR$
qui con{\"\i}ncide avec la distance $\distX$ sur $X$.

Notons que pour chaque $x \in X$ la fonction $r \mapsto [x,
\tfrac{1}{2} r]$ induit une isom{\'e}trie entre $I \subset \R$ et son
image dans $\arbXR$.  Par cons{\'e}quent, lorsque $\cS = [x, r]$ et $\cS'
= [x', r']$ sont des points de $\arbXR$ tels que $\cS \preccurlyeq
\cS'$, alors on~a $\cS' = [x, r']$~et
$$
\distAX(\cS, \cS') = \tfrac{1}{2}(\diam(\cS') - \diam(\cS)) ~,
$$
et la fonction $\tfrac{1}{2} \diam$ induit une isom{\'e}trie entre
$[\cS, \cS']$ et l'intervalle $[\tfrac{1}{2}\diam(\cS),
\tfrac{1}{2}\diam(\cS')]$ de $\R$.

En g{\'e}n{\'e}ral, pour chaque $\cS, \cS' \in \arbXR$ l'ensemble $[\cS, \cS
\vee \cS']$ (resp. $[\cS', \cS \vee \cS']$) est isom{\'e}trique {\`a} un
intervalle de $\R$ de longueur
$$
\distAX(\cS, \cS \vee \cS') = \tfrac{1}{2} (\sup \{ \cS, \cS' \} - \diam(\cS))
$$
$$
\text{(resp. }
\distAX(\cS', \cS \vee \cS') = \tfrac{1}{2} (\sup \{ \cS, \cS' \} - \diam(\cS')) 
\text{)}~,
$$
et l'ensemble $[\cS, \cS'] = [\cS, \cS \vee \cS'] \cup [\cS', \cS \vee \cS']$ est isom{\'e}trique {\`a} un intervalle de $\R$ de longueur
$$
\distAX(\cS, \cS \vee \cS') + \distAX(\cS', \cS \vee \cS')
=
\distAX (\cS, \cS') ~.
$$

On a montr{\'e} ainsi que pour chaque paire de points $\cS$ et $\cS'$ de $\arbXR$, l'ensemble $[\cS, \cS']$ est un arc topologique joignant $\cS$ et $\cS'$ et que $[\cS, \cS']$ est isom{\'e}trique {\`a} un intervalle de $\R$ de longueur $\distAX(\cS, \cS')$.
\subsection{Compl{\'e}tion de $\arbXR$}\label{ss:completion de arbXR}
On d{\'e}signe par $\arbX$ la compl{\'e}tion de $\arbXR$ par rapport {\`a} la distance $\distAX$ et on d{\'e}signe aussi par $\distAX$ l'extension de $\distAX$ {\`a} $\arbX$.

L'in{\'e}galit{\'e}
$$
|\diam(\cS) - \diam(\cS')| \le 2\distAX(\cS, \cS') ~,
$$
valable pour tout $\cS$ et $\cS'$ dans $\arbXR$, implique que la fonction $\diam$ s'{\'e}tend continuement {\`a} $\arbX$.
On d{\'e}signe cette extension aussi par  $\diam$.

\begin{lemm}\label{l:fonction diam}
La fonction $\diam : \arbX \to \R$ s'annule pr{\'e}cisement sur $X$.
\end{lemm}
\begin{proof}
  La fonction $\diam$ s'annule sur $X$.  D'autre part, soit $\cS \in
  \arbX$ tel que $\diam(\cS) = 0$ et soit $\{ [x_n, r_n] \}_{n \ge 0}$
  une suite de Cauchy dans $\arbXR$ convergent vers $\cS$.  On a donc
  $r_n = \diam(\cS_n) \to 0$ lorsque $n \to \infty$ et par cons{\'e}quent
  $\{ x_n \}_{n \ge 0}$ converge vers $\cS$ lorsque $n \to \infty$.
  Comme l'inclusion de $X$ dans $\arbXR$ est isom{\'e}trique et comme $X$
  est complet, on conclut que $\cS \in X$.
\end{proof}
\subsection{Extension des fonctions  $\cdot \vee \cdot$  et  $\sup \{ \cdot, \cdot \}$}\label{ss:extension de ^ et sup}
Le lemme ci-dessous implique que les fonctions  $\cdot \vee \cdot$  et  $\sup \{ \cdot, \cdot \} = \diam(\cdot \vee \cdot)$  s'{\'e}tendent de fa{\c c}on continue {\`a} $\arbX \times \arbX$.
Les propri{\'e}t{\'e}s et les notations d{\'e}crites aux paragraphes pr{\'e}c{\'e}dents s'{\'e}tendent par continuit{\'e} {\`a} $\arbX$.

Notons que pour $\cS, \cS' \in \arbXR$ l'identit{\'e}
$$
2\distAX(\cS, \cS \vee \cS') = \sup\{ \cS, \cS' \} - \diam(\cS)
$$
implique qu'on a $\diam(\cS) \le \sup \{ \cS, \cS' \}$, avec {\'e}galit{\'e} si et seulement si $\cS = \cS \vee \cS'$.

\begin{lemm}\label{l:extension de ^}
Pour $\cS, \cS_0, \cS_1 \in \arbXR$ on~a
$$
\distAX(\cS \vee \cS_0, \cS \vee \cS_1) \le \distAX(\cS_0, \cS_1) ~.
$$
\end{lemm}
\begin{proof}
Posons $\cS = [x, r]$ et $\cS_i = [x_i, r_i]$.
Alors,
\begin{eqnarray*}
2\distAX(\cS \vee \cS_0, \cS \vee \cS_1)
& = &
2 \max \{ r, r_0, r_1, \distX(x, x_0), \distX(x, x_1) \}
\\ & & \
- \max \{ r, r_0, \distX(x, x_0) \}
- \max \{ r, r_1, \distX(x, x_1) \}
\\ & \le &
\left|  \max \{ r_0, \distX(x, x_0) \} - \max \{ r_1, \distX(x, x_1) \}  \right| ~.
\end{eqnarray*}

Sans perte de g{\'e}n{\'e}ralit{\'e} on suppose qu'on a
$$
\max\{r_0, \distX(x, x_0) \} \ge \max \{ r_1, \distX(x, x_1) \}.
$$
Dans le cas $r_0 \ge \distX(x, x_0)$, on~a
$$
\max \{ r_0, \distX(x, x_0) \} - \max \{ r_1, \distX(x, x_1) \}
\le r_0 - r_1 \le 2\distAX(\cS_0, \cS_1) ~.
$$
On se ram{\`e}ne donc au cas $\distX(x, x_0) > r_0$~; on~a alors
$$
\distX(x, x_0) = \max\{r_0, \distX(x, x_0) \} \ge \max \{ r_1,
\distX(x, x_1) \} ~.
$$
Si $\distX(x, x_1) = \distX(x, x_0)$ alors il n'y a rien \`a
montrer~; on suppose donc que $\distX(x, x_1) < \distX(x, x_0)$.
Alors $\distX(x_0, x_1) = \distX(x, x_0) \ge \max\{ r_0, r_1 \}$ d'o{\`u},
\begin{multline*}
2\distAX(\cS_0, \cS_1)
=
2\distX(x_0, x_1) - r_0 - r_1
\\ \ge
\distX(x, x_0) - \max \{ r_1, \distX(x, x_1) \} 
= 
2\distAX(\cS \vee \cS_0, \cS \vee \cS_1) ~.
\end{multline*}
\end{proof}
\subsection{Structure d'arbre}
Rappelons qu'un espace m{\'e}trique $(\sA, d)$ est un \textit{arbre r{\'e}el}
si pour chaque paire de points $a$ et $a'$ dans $\sA$ il existe un et
un seul arc topologique joignant $a$ et $a'$ et que cet arc
topologique est isom{\'e}trique {\`a} un intervalle de $\R$ de longuer $d(a,
a')$.

\begin{prop}
  Chacun des espaces m{\'e}triques $(\arbXR, \distAX)$ et $(\arbX,
  \distAX)$ est un arbre r{\'e}el.
\end{prop}

La preuve de cette proposition est ci-dessous.  On utilisera la
notation usuelle d'intervalles pour d{\'e}signer les arcs topologiques
dans $\arbX$, de la m{\^e}me fa{\c c}on comme on~a fait pour $\arbXR$,
voir~\S~\ref{ss:segments}.  Une des cons{\'e}quences de cette proposition
est qu'une partie de $\arbXR$ (resp. $\arbX$) est connexe si et
seulement si elle est convexe.  En particulier, toute intersection
d'ensembles connexes est connexe.

La preuve de la proposition s'appuit sur le lemme suivant.

\begin{lemm}
  Soit $\cS_0 = [a_0, r_0]$ un point de $\arbXR$ et soit $\gamma : [0,
  1] \to \arbXR$ une fonction continue telle que $\gamma(0)
  \preccurlyeq \cS_0$ et telle que pour chaque $t \in [0, 1]$ on ait
  $\diam(\gamma(t)) \le r_0$.  Alors pour tout $t \in [0, 1]$ on~a
  $\gamma(t) \preccurlyeq \cS_0$ avec {\'e}galit{\'e} si et seulement si
  $\diam(\gamma(t)) = r_0$.
\end{lemm}
\begin{proof}
  Etant donn{\'e} $a \in X$ et $r > 0$ soit
$$
\arbboule^-(a, r) = \{ [a', r'] \in \arbXR \mid \max \{ \distX(a, a'), r' \} < r \} ~.
$$
Comme $\max \{ \distX(a, a'), r' \} = \sup \{ [a', r'], a \}$ et la fonction $\sup \{ \cdot, \cdot \}$ est continue, l'ensemble $\arbboule^-(a, r)$ est ouvert dans $\arbXR$.

Etant donn{\'e} $r \ge r_0$ et $a, a' \in X$, les ensembles
$\arbboule^-(a,r), \arbboule^-(a', r)$ sont soit {\'e}gaux, soit
disjoints.  Comme pour tout $r > r_0$ le point $\cS_0$ appartient {\`a}
$\arbboule^-(a_0, r)$, on conclut que $\gamma([0, 1]) \subset
\arbboule^-(a_0, r)$ et que
$$
\gamma([0, 1]) \subset \bigcap_{r > r_0} \arbboule^-(a, r) = \{ \cS \in \arbXR \mid \cS \preccurlyeq \cS_0 \} ~.
$$
Pour chaque $t \in [0, 1]$ on~a donc $2\distAX(\gamma(t), \cS_0) =
r_0 - \diam(\gamma(t))$ et par cons{\'e}quent on~a $\gamma(t) = \cS_0$ si
et seulement si $\diam(\gamma(t)) = r_0$.
\end{proof}

\

\begin{proof}[D{\'e}monstration de la proposition]
  Comme la compl{\'e}tion d'un arbre r{\'e}el est aussi un arbre r{\'e}el, il
  suffit de montrer que $(\arbXR, \distAX)$ est un arbre r{\'e}el.  Pour
  cela, il reste {\`a} montrer que pour chaque paire de points $\cS = [x,
  r], \cS' = [x', r'] \in \arbXR$, l'ensemble $[\cS, \cS']$ est le
  seul arc topologique joignant $\cS$ et $\cS'$.  Soit $\gamma : [0,
  1] \to \arbXR$ une fonction continue et injective telle que
  $\gamma(0) = \cS$ et $\gamma(1) = \cS'$.

On montrera d'abord que
$$
\rho
=
\sup \{ \diam(\gamma(t)) \mid t \in [0, 1] \} \ge \sup \{ \cS, \cS' \} ~.
$$
En effet, le lemme pr{\'e}c{\'e}dent implique que pour chaque $t \in [0, 1]$ on~a $\gamma(t) \preccurlyeq [x, \rho]$.
On a donc $\cS' \preccurlyeq [x, \rho]$ et par cons{\'e}quent $\rho \ge \max \{ r, r', \distX(x, x') \} = \sup \{ \cS, \cS' \}$.

Il reste \`a montrer que $\gamma([0, 1]) = [\cS, \cS']$.  Comme
$\gamma$ est injective, il suffit de montrer que pour tout $[x_0, r_0]
\in [\cS, \cS']$ il existe $t_0 \in [0, 1]$ tel que $\gamma(t_0) =
[x_0, r_0]$.  On se ram{\`e}ne au cas o{\`u} $[x_0, r_0] \in [\cS, \cS \vee
\cS']$~; c'est {\`a} dire qu'on peut supposer $x_0 = x$ et qu'on a $r_0
\in [r, \sup \{ \cS, \cS' \} ]$.  Alors l'infimum
$$
t_0 = \inf \{ t \in [0, 1] \mid \diam(\gamma(t)) \ge r_0 \},
$$
est bien d{\'e}fini et le lemme pr{\'e}c{\'e}dent implique qu'on~a $\gamma(t_0) = [x_0, r_0]$.
\end{proof}
\subsection{L'ordre partiel  $\preccurlyeq$  sur $\arbX$}
On prend la propri{\'e}t{\'e} ~\eqref{e-definition ordre}, valable sur
$\arbXR$, pour {\'e}tendre l'ordre partiel $\preccurlyeq$ {\`a} $\arbX$ ~:
pour deux points $\cS$ et $\cS'$ de $\arbX$ on d{\'e}finit $\cS
\preccurlyeq \cS'$ si et seulement si $\cS \vee \cS' = \cS'$.  La
transitivit{\'e} de $\preccurlyeq$ est donn{\'e}e par la partie~$2$ du lemme
suivant.

Notons que $\cS \preccurlyeq \cS'$ implique $\distAX(\cS, \cS') = \tfrac{1}{2} (\diam(\cS') - \diam(\cS))$ et par cons{\'e}quent dans ce cas on~a $\diam(\cS) \le \diam(\cS')$ avec {\'e}galit{\'e} si et seulement si $\cS = \cS'$.

\begin{lemm}\label{l:transitivite}
Soient $\cS, \cS' \in \arbX$ tels que $\cS \preccurlyeq \cS'$.
Alors on~a les propri{\'e}t{\'e}s suivantes.
\begin{enumerate}
\item[1.]
Pour $\tcS \in \arbX$ tel que $\sup \{ \cS, \tcS \} \le \diam(\cS')$ on~a $\tcS \preccurlyeq \cS'$.
\item[2.]
Pour $\cS'' \in \arbX$ tel que $\cS'' \succcurlyeq \cS'$, on~a $\cS'' \succcurlyeq \cS'$.
\end{enumerate}
\end{lemm}
\begin{proof}
\

\partn{1}
On~a
$$
\diam( \cS' )
\le
\sup \{ \tcS, \cS' \}
\le
\max \{ \sup \{ \tcS, \cS \}, \sup \{ \cS, \cS' \} \}
\le
\diam ( \cS' )~.
$$
C'est {\`a} dire qu'on~a $\distAX(\cS', \tcS \vee \cS') = \sup \{ \tcS, \cS' \} - \diam(\cS') = 0$ et donc $\tcS \vee \cS' = \cS'$.
Par cons{\'e}quent on~a $\tcS \preccurlyeq \cS'$ par d{\'e}finition de~$\preccurlyeq$.

\partn{2}
Notons qu'on~a
\begin{multline*}
\sup \{ \cS, \cS' \} = \diam( \cS \vee \cS' ) = \diam (\cS')
\le \sup \{\cS', \cS'' \} = \diam (\cS'')~.
\end{multline*}
Comme par hypoth{\`e}se $\cS'' \succcurlyeq \cS'$, la partie~$1$ implique qu'on a $\cS'' \succcurlyeq \cS$.
\end{proof}

\subsection{Propri{\'e}t{\'e}s de l'ordre partiel $\preccurlyeq$}\label{ss:proprietes ordre}
L'objectif de cette section est de montrer la proposition suivante.
\begin{prop}\label{p:proprietes ordre}
Pour $\cS$ et $\cS'$ dans $\arbX$ quelconques on~a les propri{\'e}t{\'e}s suivantes.
\begin{enumerate}
\item[1.]
Le point $\cS \vee \cS'$ est le plus petit point de $\arbX$ qui est {\`a} la fois plus grand que $\cS$ et que $\cS'$.
\item[2.]
Le segment g{\'e}od{\'e}sique $[\cS, \cS']$ de $\arbX$ joignant $\cS$ {\`a} $\cS'$ est {\'e}gal {\`a} l'ensemble
$$
\{ \tcS \in \arbX \mid \cS \preccurlyeq \tcS \preccurlyeq \cS \vee \cS'
\ \text{ ou } \
\cS' \preccurlyeq \tcS \preccurlyeq \cS \vee \cS' \} ~.
$$
\item[3.] Si $\cS' \succcurlyeq \cS$ et $\cS \neq \cS'$, alors $\cS' \in \arbXR$.
\item[4.] Si les points $\cS$ et $\cS'$ sont distincts, alors on~a $\cS \vee \cS' \in \arbXR$.
\end{enumerate}
\end{prop}
La d{\'e}monstration de cette proposition s'appui sur le lemme suivante.
\begin{lemm}\label{l:segments de infty}
Pour un point $\cS$ de $\arbX$ on~a les propri{\'e}t{\'e}s suivantes.
\begin{enumerate}
\item[1.]
Pour chaque $r \in I$ satisfaisant $r \ge \diam(\cS)$ il existe un unique point $\cS(r)$ de $\arbX$ tel que $\cS(r) \succcurlyeq \cS$ et $\diam(\cS(r)) = r$.
Lorsque $r > \diam(\cS)$ on~a $\cS(r) \in \arbXR$.
\item[2.]
L'ensemble $\{ \cS(r) \in \arbX \mid r \in I, r \ge \diam(\cS) \}$ est compl{\`e}tement ordonn{\'e} par rapport {\`a} l'ordre partiel  $\preccurlyeq$  et la fonction  $\tfrac{1}{2}\diam$  induit une isom{\'e}trie croissante entre cet ensemble et l'intervalle $\{ r \in I \mid r \ge \tfrac{1}{2}\diam(\cS) \}$ de $\R$.
\end{enumerate}
\end{lemm}
\begin{proof}
\

\partn{1}
Pour montrer l'unicit{\'e}, soient $\cS'$ et $\cS''$ des points de $\arbX$ tels que $\cS' \succcurlyeq \cS$, $\cS'' \succcurlyeq \cS$ et $r \= \diam (\cS') = \diam (\cS'')$.
Alors
$$
r \le \sup \{ \cS', \cS'' \} \le \max \{ \sup \{ \cS', \cS \}, \sup \{ \cS, \cS'' \} \} = r ~,
$$
d'o{\`u} $\cS' = \cS' \vee \cS'' = \cS''$.

Montrons maintenant l'existence de $\cS(r)$.  Lorsque $r = \diam(\cS)$
le point point $\cS$ satisfait les propri{\'e}t{\'e}s d{\'e}sir{\'e}es.  On se ram{\`e}ne
alors au cas $r > \diam (\cS)$.  Alors il existe un point $[x', r']$
de $\arbXR$ tel que $\sup \{ \cS, [x', r'] \} \le r$.  On a donc $r'
\le r$ et la partie~$1$ du Lemme~\ref{l:transitivite} implique qu'on~a
$\cS \preccurlyeq [x', r]$.  On~a montr\'e donc que le point $[x', r]
\in \arbXR$ satisfait les propri{\'e}t{\'e}s d{\'e}sir{\'e}es.

Il s'en suit de la preuve que lorsque $r > 0$, on~a $\cS(r) \in \arbXR$.

\partn{2}
Etant donn{\'e}s $r, r'$ tels que $r' \ge r \ge \diam(\cS)$, on~a
$$
\sup \{ \cS(r), \cS(r') \}
\le
\max \{ \sup \{ \cS(r), \cS \}, \sup \{ \cS, \cS(r') \} \}
= r' ~,
$$
et alors la partie~$1$ du Lemme~\ref{l:transitivite} implique qu'on~a $\cS(r) \preccurlyeq \cS(r')$.
La derni{\`e}re assertion de la proposition suit du fait que pour chaque paire de points $\tcS$ et $\tcS'$ de $\arbX$ tels que $\tcS \preccurlyeq \tcS'$, on~a $\distAX(\tcS, \tcS') = \tfrac{1}{2} (\diam(\tcS') - \diam(\tcS))$.
\end{proof}
\begin{proof}[D{\'e}monstration de la Proposition~\ref{p:proprietes ordre}]
\

\partn{1}
Lorsque $\cS$ et $\cS'$ sont des points de $\arbXR$ il est facile de voir qu'on a $\cS \vee (\cS \vee \cS') = \cS \vee \cS'$.
Par continuit{\'e} cette identit{\'e} est aussi valable lorsque $\cS$ et $\cS'$ sont des points de $\arbX$.
On a donc $\cS \preccurlyeq \cS \vee \cS'$ et par sym{\'e}trie $\cS' \preccurlyeq \cS \vee \cS'$.
Pour montrer que $\cS \vee \cS'$ est le plus petit point satisfaisant ces propri{\'e}t{\'e}s, par le Lemme~\ref{l:segments de infty} il suffit de montrer que si $\tcS$ est un point tel que $\cS \preccurlyeq \tcS$ et $\cS' \preccurlyeq \tcS$, alors $\diam(\tcS) \ge \sup \{ \cS, \cS' \}$.
En effet, dans ce cas on~a
$$
\sup \{ \cS, \cS' \} \le \max \{ \sup \{ \cS, \tcS \}, \sup \{ \tcS, \cS' \} \} = \diam(\tcS) ~.
$$

\partn{2}
Le Lemme~\ref{l:segments de infty} implique que l'ensemble $\{ \cS \preccurlyeq \tcS \preccurlyeq \cS \vee \cS' \}$ (resp. $\{ \cS' \preccurlyeq \tcS \preccurlyeq \cS \vee \cS' \}$) est un arc topologique joignant $\cS$ et $\cS \vee \cS'$ (resp. $\cS'$ et $\cS \vee \cS'$).
D'autre part, la partie~$1$ implique que l'intersection de ces arcs consiste du point $\cS \vee \cS'$.
Leur r{\'e}union est donc le seul arc topologique dans l'arbre r{\'e}el $\arbX$ joignant $\cS$ et $\cS'$.

\partn{3}
Cette assertion est une cons{\'e}quence imm{\'e}diate de la partie~$1$ du Lemme~\ref{l:segments de infty}.

\partn{4}
On a montre dans la partie~$2$ qu'on a $\cS \preccurlyeq \cS \vee \cS'$ et $\cS' \preccurlyeq \cS \vee \cS'$.
Comme les points $\cS$ et $\cS'$ sont distincts, on~a, soit $\cS \neq \cS \vee \cS'$, soit $\cS' \neq \cS \vee \cS'$.
Alors la partie~$3$ implique qu'on~a $\cS \vee \cS' \in \arbXR$.
\end{proof}

\subsection{Boules de $\arbX$}\label{ss:boules de arbX}
Etant donn{\'e} $(x, r) \in X \times I$, l'ensemble
\begin{equation}\label{e-arbboule fermee}
\{ \cS \in \arbK \mid \sup \{ x, \cS \} \le r \}
\end{equation}
ne d{\'e}pend que de la classe d'{\'e}quivalence de $(x, r)$.  On notera cet
ensemble par $\arbboule_{[x, r]}$.  Lorsque $r > 0$ on~appelera cet
ensemble \textit{boule ferm{\'e}e} de $\arbX$ et l'ensemble
\begin{equation}\label{e-arbboule ouverte}
\{ \cS \in \arbX \mid \sup \{ x, \cS \} < r \}
\end{equation}
sera appel{\'e} \textit{boule ouverte} de $\arbXR$.  De plus, on dira que le
point $[x, r]$ de $\arbXR \setminus X$ \textit{est le point d\'etermin{\'e}
  par la boule} d{\'e}finie par~\eqref{e-arbboule fermee}
ou~\eqref{e-arbboule ouverte}.  Notons que l'application $\cS \mapsto
\arbboule_\cS$ induit une bijection entre $\arbXR \setminus X$ et la
collection des boules ferm{\'e}es de $\arbX$.  D'autre part, notons que le
point de $\arbKR \setminus X$ d\'etermin{\'e} par une boule ouverte
$\arbboule$ de $\arbK$ n'appartient pas {\`a} $\arbboule$, et qu'en g{\'e}n{\'e}ral deux boules ouvertes de $\arbX$ peuvent d\'eterminer le m{\^e}me point de $\arbXR \setminus X$.

Comme la fonction $\sup\{ \cdot, \cdot \}$ est continue sur $\arbX$
(cf.~\S~\ref{ss:extension de ^ et sup}), toute boule ouverte (resp.
ferm{\'e}e) de~$\arbX$ est une partie ouverte (resp. ferm{\'e}e) de $\arbX$.
Notons d'autre part que l'intersection d'une boule de $\arbX$ avec $X$
est une boule de~$X$, et que pour chaque $\cS \in \arbXR$ on~a
$\arbboule_\cS \cap X = B_\cS$.

Soit $\arbboule$ l'ensemble d{\'e}fini par~\eqref{e-arbboule ouverte}
(resp. \eqref{e-arbboule fermee}).  Alors
l'in{\'e}galit{\'e}~\eqref{ultrametrique} implique que pour tout $\cS' \in
\arbboule$ on~a
$$
\arbboule = \{ \cS \in \arbX \mid \sup \{ \cS, \cS' \} < r \}
\text{ (resp. $\arbboule = \{ \cS \in \arbX \mid \sup \{ \cS, \cS' \} \le r \}$).}
$$
Par cons{\'e}quent, si deux boules de $\arbX$ s'intersectent, alors l'une est contenue dans l'autre. 

\begin{lemm}
Chaque boule de $\arbX$ est connexe.
\end{lemm}
\begin{proof}
Soit $\arbboule$ une boule ouverte ou ferm{\'e}e de $\arbX$ et soient $\cS$ et $\cS'$ des points dans $\arbboule$.
Il suffit de montrer que pour chaque point $\tcS$ dans $[\cS, \cS']$ on~a $\sup \{ \cS, \tcS \} \le \sup \{ \cS, \cS' \}$ ou $\sup \{ \cS', \tcS \} \le \sup \{ \cS, \cS' \}$.
On se ram{\`e}ne au cas o{\`u} $\cS \preccurlyeq \tcS \preccurlyeq \cS \vee \cS'$, cf. partie~$2$ de la Proposition~\ref{p:proprietes ordre}.
Dans ce cas on~a $\sup \{ \cS, \tcS \} = \diam(\tcS) \le \sup \{ \cS, \cS' \}$.
\end{proof}

\newpage
\section{L'arbre associ{\'e} {\`a} un corps ultram{\'e}trique}\label{s:l'arbre d'un corps}
Soit $K$ un corps muni d'une norme ultram{\'e}trique $| \cdot |$, pour
laquelle $K$ est complet et {\`a} groupe des valeurs non discret.  Dans
cette section on consid{\`e}re l'arbre r{\'e}el $\arbK$ construit dans la
section pr{\'e}c{\'e}dente, lorsque l'espace ultram{\'e}trique $(X, \distX)$ est
{\'e}gal {\`a} $K$ munit de la distance induite par la norme~$| \cdot |$.  On
consid{\`e}re aussi l'arbre projectivis{\'e} $\arbpK \= \arbK \sqcup \{ \infty
\}$ muni de la distance chordale.
\subsection{L'arbre $\arbK$}\label{ss:l'arbre de K}
Soient $I, \arbKR, \diam, \ldots $ comme dans la
Seccion~\ref{s:l'arbre d'un espace}, lorsque l'espace ultram{\'e}trique
$(X, \distX)$ est {\'e}gal au corps $K$ muni de la distance induite par
la norme $| \cdot |$.  Notons en particulier que~$I$ est l'intervalle
$[0, + \infty)$ de $\R$.  Comme par hypoth{\`e}se le groupe des valeurs de
$K$ est non discret, il est facile de voir que l'appliaction $\cS
\mapsto B_\cS$ induit une bijection entre $\arbKR \setminus K$ et la
collection des boules ferm{\'e}es et irrationnelles de~$K$.  Notons
d'autre part que l'intersection d'une boule ouverte (resp. ferm{\'e}e) de
$\arbK$ avec $K$ est une boule ouverte ou irrationnelle (resp. ferm{\'e}e
ou irrationnelle) de $K$ et que pour chaque boule ouverte ou
irrationnelle (resp. ferm{\'e}e ou irrationnelle) $B$ de~$K$ il existe une
unique boule ouverte (resp. ferm{\'e}e) $\arbboule$ de~$\arbK$ telle que
$\arbboule \cap K = B$.
\subsection{Fonction  $| \cdot |$}\label{ss:| |}
Pour $[x, r] \in \arbKR$ on pose $\left| [x, r] \right| \= \max\{ |x|,
r \}$.  La fonction $| \cdot | : \arbKR \to I$ ainsi d{\'e}finie co{\"\i}ncide
avec la norme $| \cdot |$ sur $K$.  L'in{\'e}galit{\'e}
$$
\left|  |\cS| - |\cS'|  \right| \le 2\distAX(\cS, \cS'),
$$
valable pour tout $\cS$ et $\cS'$ dans $\arbKR$, implique que la
fonction $| \cdot |$ s'{\'e}tend contin\^ument {\`a} $\arbK$.  On d{\'e}signe aussi
par $| \cdot |$ cette extension.

Notons que pour $\cS \in \arbK$ on~a $\diam(\cS) \le | \cS |$, avec
{\'e}galit{\'e} si et seulement si $\cS = [0, |\cS|]$.  Notons d'autre part
qu'on a $\cS \vee 0 = [0, |\cS|]$ et $\sup \{ \cS, 0 \} = | \cS |$.
De plus, si $B = \{ z \in K \mid |z| < r \}$ (resp. $B = \{ z \in K
\mid |z| \le r \}$), alors la boule ouverte (resp. ferm{\'e}e) de $\arbKR$
correspondante est {\'e}gale {\`a}
$$
\{ \cS \in \arbK \mid |\cS| < r \}
\text{ (resp. $\{ \cS \in \arbK \mid |\cS| \le r \}$)} ~.
$$
\subsection{Distance chordale}\label{ss:distance chordale}
On d{\'e}signe par $\Delta$ la distance sur $\arbK$ qu'on obtient
lorsqu'on int{\`e}gre la densit{\'e} $\cS \mapsto \max \{1, | \cS| \}^{-1}$
par rapport {\`a} la distance $\distAX$.  Comme cette densit{\'e} est
continue, l'espace m{\'e}trique $(\arbK, \Delta)$ est un arbre r{\'e}el.

On v{\'e}rifie ais\'ement que pour $\cS, \cS' \in \arbK$ on~a
\begin{multline*}
  \Delta(\cS, \cS') =
  \frac{\sup\{ \cS, \cS' \}}{\max \{ 1, | \cS| \} \cdot \max \{ 1, |\cS'| \}} \\
  -\tfrac{1}{2} \frac{\diam(\cS)}{\max \{1, |\cS| \}^2} - \tfrac{1}{2}
  \frac{ \diam(\cS')}{\max\{1, |\cS'| \}^2} ~.
\end{multline*}
Notons en particulier que la restriction de cette distance {\`a} $K$
co{\"\i}ncide avec la restriction de la distance chordale de $\pK$ {\`a}~$K$.

Notons d'autre part que pour chaque $r \ge 1$ et chaque paire de
points $\cS, \cS'$ dans la boule $\ov{B_r} \= \{ \tcS \in \arbK \mid
|\tcS| \le r \}$, on~a
\begin{equation}\label{comparaison}
r^{-2} \delta(\cS, \cS) \le \Delta(\cS, \cS') \le \delta(\cS, \cS')
\end{equation}
(on s'appuit ici sur le fait que la boule $\ov{B_r}$ de $\arbK$ est
connexe, cf. \S~\ref{ss:boules de arbX}).  On conclut en particulier
que les distances $\Delta$ et $\distAX$ induisent la m{\^e}me topologie
sur $\arbK$.

On d{\'e}signe par $\arbpK$ la r{\'e}union disjointe de $\arbK$ et d'un point
qu'on d{\'e}signe par $\infty$.  On id{\'e}ntifie alors la partie $K \sqcup \{
\infty \}$ de $\arbpK$ {\`a} $\pK$.

\begin{prop}\label{p:distance chordale}
  La distance $\Delta$ sur $\arbK$ s'{\'e}tend en une distance sur $\arbpK
  = \arbK \sqcup \{ \infty \}$, encore not{\'e}e $\Delta$, pour laquelle
  $(\arbpK, \Delta)$ est un arbre r{\'e}el complet.  De plus, la
  restriction de $\Delta$ {\`a} $\pK$ co{\"\i}ncide avec la distance chordale
  de $\pK$.
\end{prop}
\begin{proof}
  Par~\eqref{comparaison}, toute suite de Cauchy $\{ \cS_i \}_{i \ge
    0}$ dans $(\arbK, \Delta)$ qui n'est pas convergente satisfait $|
  \cS_i | \to \infty$ lorsque $i \to \infty$.  Comme le diam{\`e}tre de
  $\arbK \setminus \ov{B_r}$ par rapport {\`a} $\Delta$ est {\'e}gal {\`a} $\max
  \{ 1, r \}^{-1}$, toutes les suites de Cauchy dans $\arbK$ qui ne
  sont pas convergentes, sont {\'e}quivalentes.  On peut alors identifier
  la compl{\'e}tion de $\arbK$ par rapport {\`a} $\Delta$ avec l'ensemble
  $\arbpK = \arbK \sqcup \{ \infty \}$.  Il s'en suit que la distance
  $\Delta$ s'{\'e}tend {\`a} $\arbpK$ et que $\arbpK$ est un arbre r{\'e}el
  complet par rapport {\`a} cette distance (on s'appuit ici sur le fait que
  la compl{\'e}tion d'un arbre r{\'e}el est aussi un arbre r{\'e}el).
  
  Comme la restriction de la distance $\Delta$ sur $\arbK$ {\`a} $K$
  co{\"\i}ncide avec la distance chordale et comme $\pK$ est compl{\'e}t par
  rapport {\`a} la distance chordale, on conclut que la restriction de la
  distance $\Delta$ sur $\arbpK = \arbK \sqcup \{ \infty \}$ {\`a} $\pK =
  K \sqcup \{ \infty \}$ co{\"\i}ncide avec la distance chordale.
\end{proof}
\begin{rema}\label{r-arbre projectif}
  Il est facile de voir que l'arbre r{\'e}el $(\arbpK, \Delta)$ est
  isom{\'e}trique {\`a} l'arbre $\arb_{(\pK, \Delta)}$, d{\'e}finie dans la
  Section~\ref{s:l'arbre d'un espace} lorsque l'espace ultram{\'e}trique
  $(X, \distX)$ est {\'e}gal {\`a} $(\pK, \Delta)$.
\end{rema}

\subsection{Extension des fonctions $\diam$, $| \cdot |$, $\cdot \vee \cdot$, $\sup \{ \cdot, \cdot \}$ et $\preccurlyeq$}\label{ss:extension a arbpK}
On notera par $\ov{\R} \= \R \cup \{ \pm \infty \}$ la droite r{\'e}elle
achev{\'e}e.  On {\'e}tend les fonctions $\diam, | \cdot | : \arbK \to \R$ en
fonctions d{\'e}finies sur $\arbpK = \arbK \cup \{ \infty \}$ et prenant
valeurs dans $\ov{\R}$, par
$$
\diam(\infty) = | \infty | = + \infty ~.
$$
D'autre part, on {\'e}tend les fonctions  $\cdot \vee \cdot : \arbK \times \arbK \to \arbK$ et $\sup \{ \cdot, \cdot \} : \arbK \times \arbK \to \R$  en fonctions d{\'e}finies sur $\arbpK \times \arbpK$ {\`a} valeurs dans $\arbpK$ et $\ov{\R}$, respectivement, par
$$
\cdot \vee \infty = \infty \vee \cdot = \infty
\ \text{ et } \
\sup \{ \cdot, \infty \} = \sup \{ \infty, \cdot \} = + \infty ~.
$$
Notons que la fonction $\diam$ est discontinue en~$\infty$ et que
la fonction $\sup \{ \cdot, \cdot \}$ est discontinue en chaque point
de la forme $(\infty, \cdot)$ ou $(\cdot, \infty)$.
\begin{lemm}
  La fonction $| \cdot |$ est continue sur $\arbpK$ et la fonction $\cdot \vee \cdot$ est continue sur $\arbpK \times \arbpK$.
\end{lemm}
\begin{proof}
  La formule $\Delta(\cS, \infty) = \max \{ 1, | \cS | \}^{-1}$,
  valable pour tout point $\cS$ dans $\arbK$, implique que la fonction
  $| \cdot |$ est continue en $\infty$ et donc sur~$\arbpK$.  La
  continuit{\'e} de $\cdot \vee \cdot$ est alors une cons{\'e}quence de
  l'identit{\'e} $|\cS \vee \cS'| = \max \{ |\cS|, |\cS'| \}$, valable
  pour tout $\cS$ et $\cS'$ dans $\arbK$.
\end{proof}

On {\'e}tend la relation $\preccurlyeq$ {\`a} un relation d{\'e}finie sur $\arbpK$
par $\cdot \preccurlyeq \infty$.  Il est facile de voir que
$\preccurlyeq$ est une relation d'ordre partiel d{\'e}fini sur $\arbpK$.
Dans le lemme suivant on d{\'e}montre que la partie~$2$ de la
Proposition~\ref{p:proprietes ordre} s'{\'e}tend {\`a}~$\arbpK$.  On v{\'e}rifie
sans probl\`eme que les parties $1$, $3$ et~$4$ de la
Proposition~\ref{p:proprietes ordre} sont encore valables lorsqu'on
remplace $\arbK$ par $\arbpK$ et $\arbKR$ par $\arbpKR \= \arbKR \cup
\{ \infty \}$.  On conclut en particulier que pour chaque paire de
points distincts $\cS, \cS' \in \arbpK$ l'ensemble $(\cS, \cS')$ est
contenu dans $\arbpKR$.
\begin{lemm}
Pour chaque paire de points $\cS$ et $\cS'$ dans $\arbpK$ on~a
$$
[\cS, \cS']
=
\{ \tcS \in \arbpK \mid \cS \preccurlyeq \tcS \preccurlyeq \cS \vee \cS'
\ \text{ ou } \
\cS' \preccurlyeq \tcS \preccurlyeq \cS \vee \cS' \} ~.
$$
En particulier on~a $[\cS, \infty] = \{ \tcS \in \arbpK \mid \cS
\preccurlyeq \tcS \}$.
\end{lemm}

\begin{proof}
Le cas o{\`u} $\cS$ et $\cS'$ appartiennent {\`a} $\arbK$ est donn{\'e} par la partie~$2$ de la Proposition~\ref{p:proprietes ordre}.
On se ram{\`e}ne alors au cas $\cS' = \infty$.
Le Lemme~\ref{l:segments de infty} implique que l'ensemble
$$
\{ \tcS \in \arbpK \mid \cS \preccurlyeq \tcS \}
=
\{ \tcS \in \arbK \mid \cS \preccurlyeq \tcS \} \cup \{ \infty \}
$$
est un arc topologique dans $\arbpK$ joignant $\cS$ et $\infty$.
Cet ensemble est donc {\'e}gal {\`a} $[\cS, \infty]$.
\end{proof}
\subsection{Boules de $\arbpK$}\label{ss:boules de arbpK}
Une \textit{boule ouverte} (resp. \textit{ferm{\'e}e}) de $\arbpK$ est soit une boule ouverte (resp. ferm{\'e}e) de $\arbK$, soit le compl{\'e}mentaire dans $\arbpK$ d'une boule ferm{\'e}e (resp. ouverte) de $\arbK$.
Toute boule ouverte (resp. ferm{\'e}e) de $\arbpK$ est une partie ouverte (resp. ferm{\'e}e) de $\arbpK$.

Rappelons que chaque boule de $\arbK$ d\'etermine un point de $\arbKR
\setminus K$ (\S~\ref{ss:boules de arbX}).  Lorsque $\arbboule$ est
une boule de $\arbpK$ qui n'est pas une boule de $\arbK$, alors
l'ensemble $\arbboule' \= \arbpK \setminus \arbboule$ est une boule de
$\arbK$ et on dit que le point de $\arbKR \setminus K$ d\'etermin{\'e} par
$\arbboule'$ \textit{est le point determin{\'e} par la boule} $\arbboule$.
De cette fa{\c c}on chaque boule de $\arbpK$ d\'etermine un point de $\arbKR
\setminus K$.  Notons que le point de $\arbKR \setminus K$ determin{\'e}
par une boule ouverte $\arbboule$ de $\arbpK$, n'appartient pas
{\`a}~$\arbboule$.

Lorsque $\arbboule$ est une boule ouverte (resp. ferm{\'e}e) de $\arbpK$, l'ensemble $\arbboule \cap \pK$ est une boule ouverte ou irrationnelle (resp. ferm{\'e}e ou irrationnelle) de~$\pK$.
Inversement, {\`a} chaque boule ouverte ou irrationnelle (resp. ferm{\'e}e ou irrationnelle) $B$ de $\pK$ correspond une unique boule ouverte (resp. ferm{\'e}e) $\arbboule$ de $\arbpK$ telle que $\arbboule \cap \pK = B$.
\begin{prop}\label{p:boules de arbpK}
\

\begin{enumerate}
\item[1.]
Chaque boule de $\arbpK$ est connexe.
\item[2.]
Pour chaque boule $\arbboule$ de $\arbpK$, l'ensemble $\arbbouleR = \arbboule \cap \arbpKR$ est dense dans~$\arbbouleR$.
\item[3.]
Soit $\arbboule$ une boule ouverte de $\arbpK$ et soit $\cS$ le point de $\arbpKR \setminus K$ determin{\'e} par $\arbboule$.
Alors pour chaque $\cS', \cS'' \in \arbboule$ il existe $\ov{\cS} \in \arbboule$ tel que
$$
[\cS', \cS) \cap [\cS'', \cS) = [\ov{\cS}, \cS).
$$
En particulier les segments $(\cS', \cS)$ et $(\cS'', \cS)$ s'intersectent.
\end{enumerate}
\end{prop}

\begin{proof}
\

\partn{1}
Soit $\arbboule$ une boule de $\arbK$.
Comme on~a d{\'e}j{\`a} montr{\'e} que les boules de $\arbK$ sont connexes (cf. \S~\ref{ss:boules de arbX}), il suffit de montrer que le compl{\'e}mentaire de $\arbboule$ dans $\arbpK$ est connexe.
Par \S~\ref{ss:segments} il suffit donc de montrer que pour tout point $\cS$ dans $\arbpK \setminus \arbboule$ et tout point $\cS'$ de $\arbpK$ tel que $\cS \preccurlyeq \cS'$, on~a $\cS' \in \arbpK \setminus \arbboule$.
Pour montrer cette derni{\`e}re assertion, soit $\cS_0$ un point dans~$\arbboule$.
Lorsque $\diam(\cS') \ge \sup \{ \cS_0, \cS \}$ on~a $\sup \{ \cS_0, \cS' \} \ge \sup \{ \cS_0, \cS \}$ et lorsque $\diam(\cS') \le \sup \{ \cS_0, \cS \}$ on~a 
$$
\sup \{ \cS_0, \cS \}
\le
\max \{ \sup \{ \cS_0, \cS' \}, \sup \{ \cS', \cS \} \} = \sup \{ \cS_0, \cS' \}~.
$$
Dans tous les cas on~a donc $\sup \{ \cS_0, \cS' \} \ge \sup \{ \cS_0, \cS \}$ et par cons{\'e}quent $\cS \not \in \arbboule$ implique $\cS' \not \in \arbboule$.
Ceci termine la preuve de la partie~$1$.

\partn{2}
Soit $\cS \in \arbboule$ quelconque et choisisons $\cS' \in \arbboule$ distinct de $\cS$.
Comme $\arbboule$ est connexe, on~a $[\cS, \cS'] \subset \arbboule$ et comme $(\cS, \cS') \subset \arbpKR$, on conclut que $(\cS, \cS') \subset \arbbouleR$ (\S~\ref{ss:extension a arbpK}).
Par cons{\'e}quent $\cS$ est contenu dans la fermeture topologique de $(\cS, \cS') \subset \arbbouleR$.
Comme $\cS$ est un point arbitraire de $\arbboule$, on conclut que $\arbbouleR$ est dense dans~$\arbboule$.

\partn{3}
Comme $\arbpK$ est un arbre r{\'e}el l'intersection $[\cS, \cS'] \cap [\cS', \cS''] \cap [\cS'', \cS]$ est r{\'e}duite {\`a} un unique point qu'on note par $\ov{\cS}$.
L'assertion d{\'e}sir{\'e}e est alors une cons{\'e}quence inm{\'e}diate du fait que $\ov{\cS} \in [\cS', \cS''] \subset \arbboule$ et du fait que $\cS \not \in \arbboule$.
\end{proof}
\newpage
\section{Topologie fine de $\berK$}\label{s:topologie fine}
Rappelons que l'application $\cS \mapsto B_\cS$ induit une bijection entre $\arbKR \setminus K$ et la collection des boules ferm{\'e}es et irrationnelles de $K$, voir~\S~\ref{ss:espace des boules}.

\begin{theoreme}\label{t:topologie fine}
Soit $\iota : \arbKR \setminus K \to \berK$ l'application dont l'image d'un point~$\cS$ est {\'e}gale {\`a} l'unique filtre de Cauchy minimal $\fF$ dans $\berK$ satisfaisant $B_\fF = B_\cS$.
Alors on~a les propri{\'e}t{\'e}s suivantes.
\begin{enumerate}
\item[1.]  L'application $\iota$ s'{\'e}tend {\`a} une bijection uniform{\'e}ment
  continue de $(\arbpK, \Delta)$ \`a l'espace uniforme $\berK$~; on note
  cette extension encore par~$\iota$.
\item[2.]  Pour chaque boule ouverte ou irrationnelle (resp. ferm{\'e}e ou
  irrationnelle) $B$ de $\pK$, l'image par $\iota$ de la boule ouverte
  (resp. ferm{\'e}e) de $\arbpK$ correspondante {\`a} $B$ est {\'e}gale {\`a} la boule
  ouverte (resp. ferm{\'e}e) de $\berK$ correspondante {\`a}~$B$.
\end{enumerate}
\end{theoreme}

On repousse la d{\'e}monstration de ce th{\'e}or{\`e}me au
paragraphe~\S~\ref{preuve theoreme topologie fine}.  On identifie
$\arbpK$ et $\berK$ par la bijection $\iota$ donn{\'e}e par le th{\'e}or{\`e}me.
De la m{\^e}me fa{\c c}on comme on~a fait pour $\arbpK$, on utilisera la
notation usuelle d'intervalles pour d{\'e}signer les arcs dans $\berK$.
On appellera \textit{topologie faible} la topologie originale
de~$\berK$ et on appellera \textit{topologie fine} la topologie sur
$\berK$ induite par celle de~$\arbpK$.
\begin{coro}
  Toute boule de $\berK$ est connexe par rapport {\`a} la topologie faible
  et par rapport {\`a} la topologie fine.  La fronti{\`e}re par rapport {\`a} la
  topologie fine (resp. faible) d'une boule de $\berK$, est {\'e}gale au
  point de $\berK$ d{\'e}termin{\'e} par cette boule.
\end{coro}
\begin{proof}
  Comme les boules de $\arbpK$ sont connexes
  (Proposition~\ref{p:boules de arbpK}), le th{\'e}or{\`e}me pr{\'e}c{\'e}dent
  implique que les boules de $\berK$ sont connexes.  En vue du
  Lemme~\ref{l:frontiere d'une boule}, la deuxi{\`e}me assertion est une
  cons{\'e}quence imm{\'e}diate du th{\'e}or{\`e}me.
\end{proof}

\begin{coro}\label{c:ramification}
Pour un point $\cS$ de $\berK$ on~a les propri{\'e}t{\'e}s suivantes.
\begin{enumerate}
\item[1.]
Lorsque $\cS$ est rationnel ou irrationnel, les composantes connexes de $\berK \setminus \{ \cS \}$ sont les boules ouvertes de la forme $\hB$, avec $B \in \pi(\cS)$.
\item[2.]
Lorsque $\cS$ est singulier ou appartient {\`a} $\pK$, l'ensemble $\berK \setminus \{ \cS \}$ est connexe.
\end{enumerate}
\end{coro}
\begin{proof}
\

\partn{1}
Par la Proposition~\ref{p:partition d'un point} la collection $\{ \hB \mid B \in \pi(\cS) \}$ forme une partition de $\berK \setminus \{ \cS \}$.
Comme pour chaque $B \in \pi(\cS)$ l'ensemble $\hB$ est ouvert et connexe, on conclut que cette partition est la partition de $\berK \setminus \{ \cS \}$ en composantes connexes.

\partn{2}
Lorsque $\cS$ est un point de $\pK$ ou un point singulier, on peut trouver une suite d{\'e}croissante de boules ouvertes $\{ B_n \}_{n \ge 0}$ de $\pK$, telle que $\bigcap_{n \ge 0} \hB_n = \{ \cS \}$ (cf. Proposition~\ref{p:points de berK}).
Alors pour chaque entier $n \ge 0$, l'ensemble $\ov{D_n} = \berK \setminus \hB_n$ est une boule ferm{\'e}e de $\berK$ et $\{ \ov{D_n} \}_{n \ge 0}$ est une suite d'ensembles connexes tel qu'on a $\berK \setminus \{ \cS \} = \bigcup_{n \ge 0} \ov{D_n}$.
On conclut que $\berK \setminus \{ \cS \}$ est connexe.
\end{proof}

\begin{coro}
La topologie faible et la topologie fine induisent la m{\^e}me topologie sur chaque segment g{\'e}od{\'e}sique de $\berK$.
Par cons{\'e}quent, pour une partie $\berV$ de $\berK$ les propri{\'e}t{\'e}s suivantes sont {\'e}quivalentes.
\begin{enumerate}
\item[1.]
$\berV$ est convexe.
\item[2.]
$\berV$ est connexe par rapport {\`a} la topologie fine.
\item[3.]
$\berV$ est connexe par rapport {\`a} la topologie faible.
\end{enumerate}
En particulier toute partie connexe de $\berK$ est connexe par arcs et toute intersection de parties connexes de $\berK$ est connexe.
\end{coro}
\begin{proof}
  Soit $\ell$ un segement g{\'e}od{\'e}sique de $\berK$.  Avec la topologie
  fine, $\ell$ est hom\'eomorphe {\`a} un intervalle de $\R$
  (cf.~\S~\ref{ss:distance chordale}).  Il suffit alors de montrer que
  pour tout point $\cS$ de $\ell$ qui n'est pas une extr\'emit{\'e} de
  $\ell$, chacune des composantes connexes de $\ell \setminus \{ \cS
  \}$ par rapport {\`a} la topologie fine est un ouvert dans $\ell$ par
  rapport {\`a} la topologie faible.
  
  Soit alors $\cS \in \ell$ un point qui ne soit pas une extr\'emit{\'e}
  de $\ell$.  Alors $\ell \setminus \{ \cS \}$ poss{\`e}de~$2$ composantes
  connexes, qu'on d{\'e}signe par $\ell'$ et $\ell''$.  Comme $\cS$ n'est
  pas une extr\'emit{\'e} de $\ell$, il s'en suit que $\cS$ est un point
  rationnel ou irrationel de $\berK$ (voir \S~\ref{ss:extension a
    arbpK}).  Le Corollaire~\ref{c:ramification} implique alors qu'il
  existe une boule $B'$ (resp. $B''$) $\in \pi(\cS)$ telle que $\ell'
  \subset \hB'$ (resp. $\ell'' \subset \hB''$).  La partie~$3$ de la
  Proposition~\ref{p:boules de arbpK} implique alors que $B' \neq
  B''$.  On conclut donc que $\ell' = \ell \cap \hB'$ et que $\ell'' =
  \ell \cap \hB''$.  Comme chacune $\hB'$ et $\hB''$ sont des parties
  ouvertes de $\berK$ par rapport {\`a} la topologie fine, on obtient
  l'assertion d{\'e}sir{\'e}e.
\end{proof}

\subsection{D{\'e}monstration du Th{\'e}or{\`e}me~\ref{t:topologie fine}}\label{preuve theoreme topologie fine}
Soit $\iota : \arbKR \setminus K \to \berK$ comme dans l'\'enonc\'e du
Th{\'e}or{\`e}me~\ref{t:topologie fine}.  On {\'e}tend $\iota$ {\`a} $\arbpKR = \arbKR
\cup \{ \infty \}$ par $\iota(\infty) = \infty$, de telle fa{\c c}on que
$\iota$ induise l'identit{\'e} sur $\pK$.

La d{\'e}monstration du Th{\'e}or{\`e}me~\ref{t:topologie fine} s'appuit sur le
lemme suivant.  Rappelons que si $B$ est une boule ferm{\'e}e ou
irrationnelle (resp. ouverte ou irrationnelle) de $\pK$, alors
$\ov{B}$ (resp. $\hB$) est une boule ferm{\'e}e (resp. ouverte) de
$\berK$, voir \S~\ref{ss:boules et affinoides de berK}.  De plus,
rappelons que pour une boule $\arbboule$ de $\arbpK$ on pose
$\arbbouleR = \arbboule \cap \arbpKR$.
\begin{lemm}\label{l:topologie fine}
Soit $\arbboule$ une boule ouverte (resp. ferm{\'e}e) de $\arbpK$ et soit $B = \arbboule \cap \pK$ la boule de $\pK$ correspondante.
Alors on~a
$$
\iota^{-1}(\hB) = \arbbouleR
\text{ (resp. } \iota^{-1}(\ov{B}) = \arbbouleR \text{)}~.
$$
\end{lemm}
\begin{proof}
On se ram{\`e}ne au cas o{\`u} $\arbboule$ ne contient pas~$\infty$.
Alors $B$ est une boule de $K$.

Supposons d'abord que la boule $\arbboule$ est ouverte.
Soient $a_0 \in K$ et $r_0 > 0$ tels que $B = \{ |z - a_0| < r_0 \}$.
Il est facile de voir que pour $[a, r] \in \arbKR$ on~a
$$
\sup \{ [a, r], a_0 \} = \max \{ |a - a_0|, r_0 \} < r
$$
si et seulement si $B_{[a, r]} = \{ z \in K \mid |z - a| \le r \} \subset B$.
Cette derni{\`e}re propri{\'e}t{\'e} est {\'e}quivalente {\`a} ce que $B$ appartient au filtre de Cauchy minimal $\iota([z, r])$ i.e. $\iota([z, r]) \in \hB$.
On a donc $\iota^{-1}(\hB) = \arbboule \cap \arbpKR$.

Supposons maintenant que la boule $\arbboule$ est ferm{\'e}e.
Soient $a_0 \in K$ et $r_0 > 0$ tels que $B = \{ |z - a_0| \le r_0 \}$.
Alors il est facile de voir que pour $[a, r] \in \arbKR$ on~a
$$
\sup \{ [z, r], a_0 \} = \max \{ |a - a_0|, r_0 \} \le r
$$
si et seulement si $B_{[a, r]} \subset B$. Cette derni{\`e}re propri{\'e}t{\'e} est {\'e}quivalente {\`a} ce que $B$ rencontre tout {\'e}l{\'e}ment de $\iota([z, r])$ i.e. $\iota([z, r]) \in \ov{B}$.
On a donc $\iota^{-1}(\ov{B}) = \arbboule \cap \arbpKR$.
\end{proof}

\begin{proof}[D{\'e}monstration du Th{\'e}or{\`e}me~\ref{t:topologie fine}]

Comme les affino{\"\i}des ouverts forment une base pour la topologie de~$\berK$ (\S~\ref{ss:boules et affinoides de berK}), le Lemme~\ref{l:topologie fine} implique que l'application $\iota : \arbpKR \to \berK$ est continue.

\case{1.1}{$\iota$ est uniform{\'e}ment continue.}
Pour montrer que l'application $\iota : \arbpKR \to \berK$ est uniform{\'e}ment continue il suffit de montrer que pour chaque recouvrement $C \in \fC$ il existe $\varepsilon > 0$ tel qu'on ait la propri{\'e}t{\'e} suivante.
Pour toute paire de points $\cS, \cS' \in \arbpKR$ tel que $\Delta(\cS, \cS') < \varepsilon$, il existe $Y \in C$ tel qu'on ait $\cS, \cS' \in \hY$.

Notons que la partie $\partial C = \bigcup_{Y \in C} \partial \hY$ de
$\berK$ est un ensemble fini de points rationnels (cf.
Lemme~\ref{frontiere d'un affinoide}).  Etant donn{\'e} $\fF \in \partial
C$ soit $Y_\fF \in C$ tel que $\fF \in \hY_\fF$.  Comme l'application
$\iota$ est continue, il existe $\varepsilon(\fF) > 0$ tel que l'image
par $\iota$ de l'ensemble
$$
\{ \cS' \in \arbpKR \mid \Delta(\iota^{-1}(\fF), \cS') \le \varepsilon(\fF) \}
$$
soit contenue dans $\hY_\fF$.
On pose $\varepsilon = \min_{\fF \in \partial C} \varepsilon(\fF)$.

Soient $\cS, \cS' \in \arbpKR$ tels que $\Delta(\cS, \cS') <
\varepsilon$ et soit $Y \in C$ tel que $\iota(\cS) \in \hY$.  On se
ramene au cas o{\`u} $\iota(\cS') \not \in \hY$.  Comme l'ensemble $[\cS,
\cS'] \subset \arbpKR$ est connexe, on conclut qu'il existe un point
$\cS_0 \in [\cS, \cS']$ tel que $\iota(\cS_0)$ appartient {\`a} $\partial
\hY$~; on pose $\fF = \iota (\cS_0)$.  Comme $\Delta (\cS_0, \cS),
\Delta(\cS_0, \cS') \le \varepsilon \le \varepsilon(\fF)$, on conclut
que $\iota(\cS), \iota(\cS') \in \hY_\fF$.  Ceci termine la
d{\'e}monstration que l'application $\iota$ est uniform{\'e}ment continue.

\case{1.2}{$\iota$ s'{\'e}tend {\`a} une bijection uniform{\'e}ment continue de $\arbpK$ {\`a} $\berK$.}
Comme $\iota : \arbpKR \to \berK$ est uniform{\'e}ment continue, elle s'{\'e}tend {\`a} une fonction uniform{\'e}ment continue de $\arbpK$ {\`a} $\berK$, voir~\cite[II, \S~$3$, Proposition~$15$]{Bou}.
On notera cette extension aussi par $\iota$.
Il reste {\`a} montrer que $\iota$ est une bijection.

Pour montrer que $\iota$ est surjective, soit $\fF$ un point singulier
de $\berK$.  Pour chaque $r > \diam(B_\fF)$, soit $B_r$ la boule ferm{\'e}e
ou irrationnelle de $K$ de diam{\`e}tre~$r$ appartenant {\`a}~$\fF$, et soit
$\fF_r$ (resp. $\cS_r$) le point de $\berK$ (resp. $\arbKR$)
correspondant.  Il est facile de voir que $\fF_r$ converge vers $\fF$
lorsque $r \to \diam(\fF)$.  D'autre part, notons que la collection
$\{ B_r \}_{r > \diam(B_\fF)}$ est compl{\`e}tement ordonn{\'e}e par rapport {\`a}
l'inclusion (Proposition~\ref{p:points de berK}) et que l'application
$\rho \mapsto \cS_{2\rho}$ est une isom{\'e}trie croissante entre $\{ \rho
\in \R \mid \rho > \tfrac{1}{2} \diam(B_\fF) \}$ et son image dans
$\arbKR$.  Par cons{\'e}quent $\cS_r$ converge vers un point $\cS$ de
$\arbK$ lorsque $r \to \diam(B_\fF)$.  Par continuit{\'e} de $\iota$ on~a
$\iota(\cS) = \fF$.

Pour montrer que $\iota$ est injective, soient $\cS, \cS' \in \arbpK$
distincts.  Pour chaque $r \in [\diam(\cS), + \infty]$ (resp. $r \in
[\diam(\cS'), + \infty]$) soit $\cS_r$ (resp. $\cS'_r$) le point de
$\arbpK$ tel que $\cS_r \succcurlyeq \cS$ (resp. $\cS_r \succcurlyeq
\cS'$) et $\diam(\cS_r) = r$.  Alors pour $r \neq \diam(\cS)$ (resp.
$r \neq \diam(\cS')$) on~a $\cS_r \in \arbpKR$ (resp. $\cS'_r \in
\arbpKR$), voir Lemme~\ref{l:segments de infty}.

Supposons d'abord $\cS$ et $\cS'$ ne sont pas reli\'es par
$\preccurlyeq$, de telle fa{\c c}on que $\cS \vee \cS'$ soit distinct de
$\cS$ et de $\cS'$.  Si l'on pose $\rho = \sup \{ \cS, \cS' \}$, alors
on~a $\rho > \max \{ \diam(\cS), \diam(\cS') \}$ et $\cS_\rho =
\cS'_\rho = \cS \vee \cS'$.  Fixons $r_0 \in (\max \{ \diam(\cS),
\diam(\cS') \}, \sup\{ \cS, \cS' \} )$ et posons
$$
\arbboule = \{ \tcS \in \arbpK \mid \sup \{ \tcS, \cS \} \le r_0 \},
B = \arbboule \cap \pK,
$$
$$
\arbboule' = \{ \tcS \in \arbpK \mid \sup \{ \tcS, \cS' \} \le r_0 \}
\text{ et }
B' = \arbboule' \cap \pK.
$$
Notons que les sous-ensembles $\ov{B}$ et $\ov{B'}$ de $\berK$ sont disjoints.
D'autre part, pour tout $r \in (\diam(\cS), r_0]$ (resp. $r \in (\diam(\cS'), r_0]$) on~a $\cS_r \in \arbboule \cap \arbpKR$ (resp. $\cS_r' \in \arbboule \cap \arbpKR$) et par cons\'equent on~a $\iota(\cS_r) \in \ov{B}$ (resp. $\cS_r' \in \ov{B'}$) (Lemme~\ref{l:topologie fine}).
Par continuit{\'e} de $\iota$ on~a alors $\iota(\cS) \in \ov{B}$ et $\iota(\cS') \in \ov{B'}$, d'o{\`u} $\iota(\cS) \neq \iota(\cS')$.

Supposons maintenant que $\cS$ et $\cS'$ sont reli\'es par $\preccurlyeq$.
On se ram{\`e}ne au cas o{\`u} $\cS \preccurlyeq \cS'$.
On a alors $\diam(\cS) < \diam(\cS')$ et $\cS' = \cS_{\diam(\cS')}$.
De plus le Lemme~\ref{l:segments de infty} implique que $\cS' \in \arbpKR$.
Fixons $r_0 \in (\diam(\cS), \diam(\cS'))$ et posons
$$
\arbboule = \{ \tcS \in \arbpK \mid \sup \{ \tcS, \cS \} \le r_0 \}
\text{ et }
B = \arbboule \cap \pK.
$$
Comme pour tout $r \ \in (\diam(\cS), r_0]$ on~a $\cS_r \in \arbboule \cap \arbpKR$, le Lemme~\ref{l:topologie fine} implique qu'on a $\iota(\cS_r) \in \ov{B}$.
Par continuit{\'e} de $\iota$ on~a $\iota(\cS) \in \ov{B}$.
D'autre part, comme $\cS' \not \in \arbboule$ le Lemme~\ref{l:topologie fine} implique qu'on a $\iota(\cS') \not \in \ov{B}$.
On a donc $\iota(\cS) \neq \iota(\cS')$.
Ceci termine la d{\'e}monstration du fait que $\iota$ est une bijection.

\partn{2} Etant donn{\'e}e une boule $\arbboule$ de $\arbpK$, soit $B =
\arbboule \cap \pK$, $\arbboule' = \arbpK \setminus \arbboule$ et $B'
= \arbboule' \cap \pK$.  Comme l'ensemble $\arbbouleR$ est dense dans
$\arbboule$ (partie~$2$ de la Proposition~\ref{ss:boules de arbpK}),
le Lemme~\ref{l:topologie fine} et la continuit{\'e} de $\iota$ impliquent
que $\iota(\arbboule)$ est contenu dans $\ov{B}$.  Par le m{\^e}me
raisonnement on~a $\iota(\arbboule') \subset \ov{B'}$.  Quitte {\`a}
remplacer $\arbboule$ par $\arbboule'$, on se ram{\`e}ne au cas o{\`u} la
boule $\arbboule$ de $\arbpK$ est ouverte.  Soit $\fF$ le point de
$\berK$ determin{\'e} par~$B$.  On~a alors $\fF \not \in \hB$ et $\ov{B} =
\hB \cup \{ \fF \}$ (Lemme~\ref{l:frontiere d'une boule}).  Comme
l'application $\iota$ est une bijection le Lemme~\ref{l:topologie
  fine} implique qu'on~a $\iota^{-1}(\fF) \not\in \arbboule$, d'o{\`u} on
conclut que $\iota(\arbboule) \subset \hB$.  Comme les ensembles
$\arbboule$ et $\arbboule'$ sont compl\'ementaires dans $\arbpK$ et les
ensembles $\hB$ et $\ov{B'}$ sont compl{\'e}mentaires dans $\berK$, on
conclut que $\iota(\arbboule) = \hB$ et $\iota(\arbboule') = \ov{B'}$.
\end{proof}

\backmatter
\newpage
\bibliographystyle{plain}

\begin{thebibliography}{FvdP}

\bibitem[BR1]{BR1} M. Baker, R. Rumely.
\textit{Analysis and dynamics on the Berkovich projective line.}
\texttt{arxiv.org/math.NT/0407433}

\bibitem[BR2]{BR2} M. Baker, R. Rumely.
\textit{Equidistribution of small points, rational dynamics, and potential theory.}
\texttt{arxiv.org/math.NT/0407426}






\bibitem[Ben1]{Be Ahlfors} R. Benedetto.
\textit{Non-Archimedean holomorphic maps and the Ahlfors Islands theorem.}
Amer. J. Math.  {\bf 125}  (2003), 581-622.

\bibitem[Ben2]{Be Ahlfors II} R. Benedetto.
\textit{An Ahlfors Islands Theorem for non-archimedean meromorphic functions.}
\texttt{arxiv.org/math.NT/0407142}


\bibitem[Ber1]{Berrouge} V.G. Berkovich.
\textit{Spectral theory and analytic geometry over non-Archimedean fields.}
Mathematical Surveys and Monographs {\bf 33}, Amer. Math. Soc., Providence, RI, $1990$.

\bibitem[Ber2]{Berblue} V.G. Berkovich.
\textit{{\'E}tale cohomology for non-Archimedean analytic spaces.}
Inst. Hautes {\'E}tudes Sci. Publ. Math.  {\bf 78} (1993), 5-161 (1994).

\bibitem[Ber3]{BerICM} V.G. Berkovich.
\textit{$p$-adic analytic spaces}.
Proceedings of the International Congress of Mathematicians, Vol. II (Berlin, 1998).  Doc. Math.  $1998$,  Extra Vol. II, 141-151

\bibitem[Bez]{Bez} J.P. B{\'e}zivin.
\textit{Dynamique des fractions rationnelles $p$-adiques.}
\texttt{www.math.unicaen.fr/$\sim$bezivin/dealatex.pdf}

\bibitem[Bou]{Bou} N. Bourbaki.
\textit{{\'E}l{\'e}ments de math{\'e}matique. Topologie g{\'e}n{\'e}rale.}
Hermann, Paris, 1971.


\bibitem[BHM]{BHM} K. Boussaf, M. Hemdaoui, N. Ma{\"\i}netti.
\textit{Tree structure on the set of multiplicative semi-norms of Kranser algebras $H(D)$.}
Rev. Mat. Complut. {\bf 13} (2000), 85-109.



\bibitem[CL]{CL} A. Chambert-Loir.
\textit{Mesures et {\'e}quidistribution sur les espaces de Berkovich.}
\texttt{arxiv.org/math.NT/0304023}

\bibitem[C]{Cherry} W. Cherry.
\textit{Non-Archimedean analytic curves in abelian varieties.}
Math. Ann. {\bf 300} (1994), 393-404.

\bibitem[CR]{CherryRu} W. Cherry, M. Ru.
\textit{Rigid analytic Picard theorems.}
Amer. J. Math.  {\bf 126}  (2004), 873-889.

\bibitem[D]{D} A. Ducros.
\textit{Espaces analytiques $p$-adiques au sens de Berkovich.}
Expos{\'e} 958 du s{\'e}minaire Bourbaki (mars 2006).




\bibitem[Es1]{Es80} A. Escassut.
\textit{The ultrametric spectral theory.}
Period. Math. Hungar. {\bf 11} (1980), 7-60.

\bibitem[Es2]{Es95} A. Escassut.
\textit{Analytic elements in p-adic analysis.}
World Scientific Publishing, 1995.

\bibitem[Es3]{Es03} A. Escassut.
\textit{Ultrametric Banach algebras.}
World Scientific Publishing Co., Inc., River Edge, NJ, 2003.

\bibitem[FJ]{FJ} C. Favre, M. Jonsson.
\textit{The valuative tree.}
Lecture Notes in Mathematics, 1853. Springer-Verlag, Berlin, 2004.

\bibitem[FR1]{FR1} C. Favre, J. Rivera-Letelier.
\textit{Th{\'e}or{\`e}me d'{\'e}quidistribution de Brolin en dynamique $p$-adique.}
C. R. Math. Acad. Sci. Paris  {\bf 339}  (2004), 271-276.

\bibitem[FR2]{FR2} C. Favre, J. Rivera-Letelier.
\textit{Equidistribution des points de petite hauteur.}
Math. Ann. {\bf 335} (2006), 311-361.

\bibitem[FvdP]{FvdP} J. Fresnel, M. van der Put.
\textit{Rigid analytic geometry and its applications.}
Progress in Mathematics, 218. Birkh{\"a}user Boston, Inc., Boston, MA, 2004.

\bibitem[Ga]{Ga} G. Garandel. 
\textit{Les semi-normes multiplicatives sur les alg{\`e}bres d'{\'e}l{\'e}ments
analytiques au sens de Krasner.}
Indag. Math. {\bf 37} (1975), 327-341.

\bibitem[Gu]{Gu} B. Guennebaud.
\textit{Sur une notion de spectre pour les alg\`ebres norm\'ees ultram\'etriques.}
Th\`ese de doctorat, Universit\'e de Poitiers, 1973.


\bibitem[Hs]{Hs} L.C. Hsia.
\textit{Closure of periodic points over a non-Archimidean field}.
J. London Math. Soc. {\bf 62} (2000), 685-700.

\bibitem[Hub]{Hub} R. Huber.
\textit{{\'E}tale cohomology of rigid analytic varieties and adic spaces.}
Aspects of Mathematics, E30. Friedr. Vieweg \& Sohn, Braunschweig, $1996$.

\bibitem[Hug]{Hug} B. Hughes.
\textit{Trees and ultrametric spaces: a categorical equivalence.}
Advances in Mathematics {\bf 189} (2004), 148­191.

\bibitem[L-S]{L-S} B. Le Stum.
\textit{Espaces de Berkovich (d'apr\`es Tate, van~der~Put, Berkovich, Schneider et Deligne).}\\
\texttt{perso.univ-rennes1.fr/bernard.le-stum/Documents/berkovich.pdf}


\bibitem[Le]{Le} A.J. Lemin.
\textit{The category of ultrametric spaces is isomorphic to the category of complete, atomic, tree-like, and real graduated lattices $LAT\sp *$.}
Algebra Universalis  {\bf 50}  (2003), 35-49.

\bibitem[M1]{Mai1} N. Ma{\"\i}netti.
\textit{Sequential compactness of some analytic spaces.}
J. Anal.  {\bf 8}  (2000), 39-54.

\bibitem[M2]{Mai2} N. Ma{\"\i}netti.
\textit{Metrizability of some analytic affine spaces.}
$p$-adic functional analysis (Ioannina, 2000),  219--225, Lecture Notes in Pure and Appl. Math., 222, Dekker, New York, 2001. 





\bibitem[vdP]{vdP} M. van der Put.
\textit{Cohomology on affinoid spaces.}
Compositio Math. {\bf 45} (1982), 165-198.

\bibitem[PS]{PS} M. van der Put, P. Schneider.
\textit{Points and topologies in rigid geometry.}
Math. Ann. {\bf 302} (1995), 81-103.


\bibitem[R-L1]{these} J. Rivera-Letelier.
\textit{Dynamique des fonctions rationnelles sur des corps locaux.}
Geometric methods in dynamics. II.  Ast{\'e}risque  {\bf 287} (2003), 147-230.

\bibitem[R-L2]{hyp}
J. Rivera-Letelier.
\textit{Espace hyperbolique $p$-adique et dynamique des fonctions rationnelles.}
Compositio Math. {\bf 138} (2003), 199-231.

\bibitem[R-L3]{periodique}
J. Rivera-Letelier.
\textit{Points p{\'e}riodiques des fonctions rationnelles dans l'espace hyperbolique $p$-adique.}
Comment. Math. Helv. {\bf 80}  (2005), 593-629.

\bibitem[R-L4]{injective} J. Rivera-Letelier.
\textit{Une caract{\'e}risation des fonctions holomorphes injectives en analyse ultram{\'e}trique.}
C. R. Math. Acad. Sci. Paris  {\bf 335}  (2002), 441-446.

\bibitem[R-L5]{elements}
J. Rivera-Letelier.
\textit{Bi-analytic elements and partial isometries of hyperbolic space.}
Ultrametric functional analysis (Nijmegen, 2002),  319-343, Contemp. Math., 319, Amer. Math. Soc., Providence, RI, 2003.

\bibitem[R-L6]{Julia/Fatou} J. Rivera-Letelier.
\textit{Th{\'e}orie de Julia et Fatou sur la droite projective de Berkovich.}
En preparation.




\bibitem[Ta]{Ta} J. Tate.
\textit{Rigid analytic spaces.}
Invent. Math. {\bf 12}  (1971), 257-289. 

\bibitem[Th]{Th} A.~Thuillier
\textit{Th\'eorie du potentiel sur les courbes en g\'eom\'etrie analytique non archim\'edienne.}
Th{\`e}se de l'universit{\'e} de Rennes, 2005.

\bibitem[W]{W} A. Werner.
\textit{Compactification of the Bruhat-Tits building of PGL by seminorms.}
Math. Z. \textbf{248} (2004), 511-526.



\end{thebibliography}

\end{document}